\documentclass[trans]{IEEEtran}
\usepackage[cmex10]{amsmath}
\usepackage{amsfonts,mathtools}
\usepackage{algorithmic}
\usepackage[]{algorithm}
\usepackage{array}
\usepackage{url}
\usepackage{bm}
\usepackage{graphicx}
\usepackage{mathrsfs}
\usepackage{amsthm} 
\usepackage{amsmath,bm}
\usepackage{subfigure}
\usepackage{cite}
\usepackage{color}
\usepackage[font={small}]{caption}
\usepackage{epstopdf}
\usepackage{amssymb}
\usepackage{lipsum}
\usepackage[normalem]{ulem}

\usepackage{fancyhdr,lipsum}
\fancypagestyle{header}{%
	\fancyhf{} % clear all fields
	
}%

\makeatletter
\let\ps@IEEEtitlepagestyle\ps@header
\makeatother

\makeatletter
\newlength \figwidth
\if@twocolumn
\else
\fi
\makeatother

\theoremstyle{plain}% default
\newtheorem{thm}{Theorem}
\newtheorem{lem}{Lemma}

\theoremstyle{definition}

\newtheorem{exmp}{Example}
\newtheorem{assumption}{}

\theoremstyle{remark}

\providecommand{\abs}[1]{\lvert#1\rvert}												% abs operator
\renewcommand{\b}[1]{\ensuremath{\mathbf{#1}}}		 							% bold
\newcommand{\Ex}[1]{\ensuremath{\mathbb{E}[#1]}} 		 % expectation operator
\newcommand{\ind}{1\hspace{-1.6mm}1}														% indicator function
\newcommand{\norm}[1]{\ensuremath{\left\|#1\right\|}}						% norm operator
\providecommand{\tr}[1]{\text{tr}\left(#1\right)}

\newcommand{\arrowAS}[1]{\ensuremath{\stackrel{\text{a.s.}}{\longrightarrow}}}

\providecommand{\norm}[1]{\left \| #1 \right \|}

\providecommand{\inPro}[2]{\left \langle #1, \; #2 \right \rangle}

\providecommand{\norm}[1]{\left \| #1 \right \|}

\providecommand{\ip}[2]{\langle #1, #2 \rangle}

\newcommand{\px}[2]{\ensuremath{\mathcal{P}_{#1}\left(#2\right)}} 		 % projection operator

\def \Rn {{\mathbb{R}}}
\def \x {{\b{x}}}
\def \z {{\b{z}}}
\def \u {{\b{u}}}
\def \y {{\b{y}}}

\def \I {{\b{I}}}

\def \v {{\b{v}}}

\def \A {{\b{A}}}

\def \E {{\mathbb{E}}}

\def \leg {{L}}
\def \Xs {{\mathcal{X}}}
\def \N {{\mathcal{N}}}

\def \cX {{\mathcal{X}}}
\def \xib {{\boldsymbol{\xi}}}
\def \e {{\b{e}}}
\def \S {{\mathcal{S}}}

\def \g {{\tilde{g}}}
\def \f {{\tilde{f}}}
\def \P {{\mathbf{P}}}
\def \pk {\mathbf{P}_k}
\def \cG {{\mathcal{G}}}
\def \cK {{\mathcal{K}}}

\def \cE {{\mathcal{E}}}

\def \h {{\tilde{h}}}
\def \lt {{\tilde{\ell}}}
\def \ut {{\tilde{u}}}
\def \dz {{\delta\z}}
\def \dx {{\delta\x}}
\def \In {{\mathbb{I}}}
\def \ot {{\mathcal{O}\left(E_T+\frac{1}{T}\right)}}

\def \zb {{\bar{\z}}}

\begin{document}
%%%%%%%%%%%%%%%%%%%%%%%%%%%%
\title{Distributed Inexact Successive Convex Approximation ADMM: Analysis-Part I}
\author{Sandeep~Kumar, Ketan~Rajawat,~and~Daniel~P.~Palomar,~\IEEEmembership{Fellow,~IEEE} \\
\thanks{ Sandeep Kumar and Daniel~P.~Palomar are with the Hong Kong University of Science and Technology (HKUST),
Hong Kong. E-mail: \{eesandeep, palomar\}@ust.hk, and Ketan Rajawat is with the Department of Electrical Engineering, Indian Institute of Technology Kanpur, India, E-mail: ketan@iitk.ac.in.}}
\maketitle

\vspace{-0.8cm}

%\begin{center}
%{\small{\bf Submitted:} \today}
%\end{center}

%\vspace{0.3cm}
{\color{black}
%%%%%%%%%%%%%%%%%%%%%%%%%%%%%%%%%%%%%%%%%%
\begin{abstract}
In this two-part work, we propose an algorithmic framework for solving non-convex problems whose objective function is the sum of a number of smooth component functions plus a convex (possibly non-smooth) or/and smooth (possibly non-convex) regularization function. The proposed algorithm incorporates ideas from several existing approaches such as alternate direction method of multipliers (ADMM), successive convex approximation (SCA), distributed and asynchronous algorithms, and inexact gradient methods. Different from a number of existing approaches however, the proposed framework is flexible enough to incorporate a class of non-convex objective functions, allow distributed operation with and without a fusion center, and include variance reduced methods as special cases. Remarkably, the proposed algorithms are robust to uncertainties arising from random, deterministic, and adversarial sources. The part I of the paper develops two variants of the algorithm under very mild assumptions and establishes first order convergence rate guarantees. The proof developed here allows for generic errors and delays, paving the way for different variance-reduced, asynchronous, and stochastic implementations, outlined and evaluated in the part II.

%First-order gradient based methods are well-regarded approach for solving a broad range of optimization problems. In many practical situations, first-order algorithms only have access to gradients with some errors stochastic or deterministic.
%	
%Recently, there has been a lot of interest in the study of first-order inexact gradient based methods under more general sources of gradient error, stochastic or deterministic for convex optimization problems. 
%	
%Large scale non-convex optimization lies at the core of several applications in science and engineering. The surge of big data for machine learning applications and ongoing expansion of networked-system has led to significant interest in provably convergent, low cost, distributed, asynchronous, and robust algorithms for scaling computations and communications. Naturally, non-convexity along with the associated constraints of big data and networked objective complicates the design and analysis of algorithms significantly. In this work we develop a novel optimization framework at the unification of successive convex approximation, inexact gradients, and alternate direction method of multipliers (ADMM), leading to a gamut of provably convergent distributed algorithms, which can be customized 
%to solve a broad range of well researched and new problems of machine learning, communication, signal processing, control and networking problems. 
\end{abstract}

\section{Introduction}

%There is a large amount of literature dedicated to the convergence analysis of these methods under the true gradient information.

{\huge{N}}on-convex optimization problems arise in a variety of signal processing, statistics, machine learning, communications, controls, and networking applications \cite{jain2017non}. As with most modern systems involving big data, the ensuing problems are also high-dimensional and large-scale, and necessitate scalable solutions \cite{krizhevsky2012imagenet}. Parallel and distributed algorithms are regarded as an enabling tools capable of meeting the challenges of big data \cite{krizhevsky2012imagenet,scutari2018parallel}. While attractive in theory, the implementation of such algorithms in distributed systems is still marred by practical issues such as heterogeneity between processor nodes, communication delays, and inexact updates. Networked systems that ignore these issues are ultimately suboptimal or rely on empirical approximations with no theoretical basis or guarantees \cite{kumar2017asynchronous,yang2016parallel}.

%{\huge{N}}on-convex optimization problems are important in a variety of signal processing, statistics, machine learning, communication, control, and networked system applications\cite{jain2017non}. Many of the modern applications require to process big data (large scale and high dimensional data) giving rise to large scale optimization \cite{krizhevsky2012imagenet}. 
%Parallel and distributed algorithms are regarded as an enabling tool to cope up with the challenges of big data \cite{krizhevsky2012imagenet,scutari2018parallel}. However, analysis and design of distributed algorithms bring additional complexities due to various factors including lack of a centralized coordinator; uncertainties introduced due to the distributed computation; and asynchrony and delay due to the heterogeneity of the nodes, to name a few. Such challenges arise even in moderate-scale problems solved over other networked systems\cite{kumar2017asynchronous,yang2016parallel}. 

%First-order gradient methods are a well-suited approach for designing distributed algorithms. They are popular because of their cheap iterations and light dependence on the problem dimension and the data size. 

Scalable optimization algorithms for solving non-convex problems generally rely on gradients to carry out the updates \cite{allen2016variance,reddi2016stochastic,huang2016stochastic}. Recent years have witnessed the development of scalable majorization and successive convex approximation algorithms that can be viewed as generalizations of the first order algorithms capable of better exploiting the structure of the loss function \cite{sun2017majorization, scutari2014parallel}. Most of these algorithms however require exact information about the loss function or its gradient, and cannot tolerate uncertainty. 

In practice, exact information about the loss function may not be available instantaneously. In real-time embedded applications for instance, the loss function (and consequently its gradient) may depend on noisy measurements of various quantities. Random noise may also get introduced into the updates via the message-passing step in networked systems, where the communication medium may be lossy or error-prone. In the context of large-scale empirical risk minimization problems, the gradient cannot be calculated at every iteration, and instead stochastic approximations of the gradient are used, as in the stochastic gradient descent (SGD) method and its variants \cite{bertsekas1999nonlinear,rosasco2014convergence,roux2012stochastic}. Beyond these, systematic errors may arise in poorly calibrated or low-quality equipment such as inertial measurement units. Delays arising from computational or communication issues may give rise to delayed gradients, which can be viewed as erroneous versions of the current gradient. If the updates involve solving an intermediate problem, such as in proximal, mirror descent, or dual-descent methods, errors may arise from solving these inner sub-problems inaccurately, e.g., due to early termination of inner iterations \cite{devolder2014first}. Errors may be modeled as adversarial if they depend on values whose measurement incurs some delay, such as the location of a target being tracked. Finally, in some settings, noise is intentionally introduced into the updates so as to ensure certain properties, e.g., noise is deliberately added to the updates to help the iterates escape saddle points \cite{jin2017escape} and quantization or sparsification may be applied to increase the communication and computational efficiency \cite{alistarh2017qsgd,wangni2017gradient}. 

Convex optimization in the presence of uncertainty or errors has been widely studied \cite{subsgrad2018,aybat2018robust,hu2017analysis,devolder2014first,schmidt2011convergence,d2008smooth,mohammadi2019robustness}. However, the effect of uncertainty on the performance of non-convex optimization algorithms has not been fully understood. Pertinent approaches either require exact gradient information or allow specific types of uncertainties. For example, distributed alternating direction method of multipliers (ADMM) methods \cite{admm,kumar2017asynchronous} converge under asynchrony but require exact gradients. Conversely, stochastic algorithms \cite{huang2016stochastic,hajinezhad2016nestt} allow unbiased gradient noise but require centralized synchronous implementation. The Successive Convex Approximation (SCA) method in \cite{scutari2018parallel} is distributed, but cannot handle errors and asynchrony. The synchronous extensions in \cite{koppel2018parallel,yang2016parallel,liu2018stochastic} can handle random noise but not uncertainty arising from deterministic and adversarial sources. More recently, inexact gradient methods with arbitrary errors have been developed \cite{sun2017convergence,dvurechensky2017gradient} but are not amenable to distributed implementations. 

In this paper we develop a generic optimization framework called \textsf{D}istributed \textsf{I}nexact \textsf{S}uccessive \textsf{C}onvex \textsf{A}pproximation \textsf{ADMM} (DISC-ADMM). The proposed approach draws from several existing approaches, such as distributed and asynchronous ADMM, distributed SCA, and Inexact Gradient Methods (IGM). The ADMM-style updates impart flexibility to the algorithm, allowing them to be implemented in distributed networks, both with and without a fusion center. Depending on the implementation, the proposed framework allows both missed as well as delayed updates. Non-convex components in the loss function and the regularizer are handled within the rubric of SCA, whose versatility enables the proposed algorithm to better exploit hidden convexity. Finally, the updates allow arbitrary error sequences and the eventual performance is characterized as a function of the total error incurred so far, regardless of its source. The proposed approach generalizes the stochastic ADMM algorithms and can be customized for a particular application, yielding a broad range of new algorithms.

The present work focuses on developing the proposed algorithm and characterizing its performance in terms of the average error norm. The implementation details, including the development of variance-reduced and/or quantized versions, and numerical tests for various empirical risk minimization problems is deferred to the part II of this paper. 

We summarize the contributions as follows:
\begin{itemize}
	\item An ADMM algorithm for solving problems with non-convex objective functions is proposed. The proposed algorithm uses updates where the non-convexity is handled within SCA framework, thereby subsuming both first order ADMM and MM-based ADMM algorithms \cite{lu2018unified, kumar2017asynchronous}. The use of such surrogate functions within ADMM is also a novel aspect which is of independent interest and can be integrated with other non-convex optimization algorithms. 	
	\item The convex surrogate functions used in lieu of the loss function can be inexact and based on old iterates. The resulting error in the gradient of the surrogate function is allowed to be arbitrary but is shown to be critical in determining the convergence of the proposed algorithm. 
	\item We establish the iterates generated by the algorithms converge to a set of first order $\epsilon-$stationary points, where $\epsilon$ is proportional to the average error norm. The use of delayed updates in constructing the surrogate functions has no effect on the convergence, provided that these delays are bounded. 
	\item Since the proposed algorithm entails ADMM-style updates, it inherits the flexibility of ADMM and its amenability to different distributed architectures, such as fusion-centric and fully distributed. 
\end{itemize}
It is remarked that the proposed algorithm subsumes existing stochastic ADMM and its variance-reduced versions, as well as algorithms involving gradient averaging or perturbations. Unlike the distributed-SCA algorithm of \cite{scutari2014parallel} which involves nested loops requiring a distributed optimization problem to be solved within another, the proposed algorithm is single loop.}

The notation used in this paper is as follows. Bold upper(lower) case letters denote matrices (vectors). $\I$ is the $N \times N$ identity matrix, $0$ denotes the all-zero matrix or vector, and $1$	denotes the all-one matrix or vector, depending on the context. For a vector $x, \norm{\x}_2$ denotes its 2 norm. For a matrix $\A, \norm{\A}$	denotes its Frobenius norm, $\norm{\A}_2$ denotes the 2 form, $\tr\A$ the trace of the matrix.For a function with multiple arguments, $\nabla$ always denotes the gradient with respect to the first argument. Likewise, $\partial$ always denotes the subgradient with respect to the first argument. For a multi-input function such as $g(\u,\v)$, we will use the compact notation $\nabla_\u g(\x,\y):= \nabla_{\u} g(\u,\v) \mid_{\u = \x, \v = \y}$. A function $f$ is said to be $L$-smooth if $\nabla f$ is Lipschitz continuous with constant $L$. Likewise, a function $f$ is said to be $m$-convex if $f$ is strongly convex with parameter $m$. 
	%%%%%%%%%%%%%%%%%%%%%%%%%%%%%%%%
	\section{Problem Formulation}
	%%%%%%%%%%%%%%%%%%%%%%%%%%%%%
%	\begin{align}
%	\mathsf{P} =	\min_{\{\x_k\},\z\in \cX} & h(\z) + \sum_{k=1}^K g_k(\x_k) \tag{$\mathcal{P}$} \label{unif}\\
%	& \text{s. t. } \pk(\x_k - \z) = 0 \ \ \ \ k = 1, \ldots, K \tag{$\mathcal{C}$} \label{unifc}
%\end{align}	

We consider a network of $K$ nodes or agents, each associated with a non-convex function $g_k:\Rn^n\rightarrow\Rn$. The agents seek to cooperatively solve the following general non-convex consensus problem:
	\begin{align}
	%\begin{array}{ll}
	\mathsf{P} =	\underset{\{\x_k\},\z\in \cX}{\text{{minimize}}} &
	\begin{array}{c}
	\hspace{.2cm}	h(\z) + \sum_{k=1}^K g_k(\x_k) \label{unif}
	\end{array}\\
	\text{{subject to}} & \begin{array}[t]{l}
		\hspace{.2cm}	\pk(\x_k - \z) = 0 \ \ \ \ k = 1, \ldots, K \label{unifc}	
	\end{array}
	%\end{array}
\end{align}
where $\x_k\in \Rn^n$, $\z \in \cX \subset \Rn^n$, and $\pk \in \Rn^{n \times n}$ is any orthogonal projection matrix such that $\sum_{k=1}^K \pk$ is positive definite. Each function $g_k$ is smooth and non-convex and takes the form $g_k(\x_k) = g_k(\pk\x_k)$, i.e., depends on $\x_k$ only through $\pk\x_k$. The function $h$ may have both smooth non-convex and convex non-smooth components. The problem in \eqref{unif} also allows hybrid scenarios, such as $h$ being sum of smooth non-convex and non-smooth convex components, $g_k$ expressible as sums of non-convex functions, and arbitrary orthogonal projection matrices. 

\begin{figure}[H]
	\begin{centering}
		\includegraphics[scale=0.33]{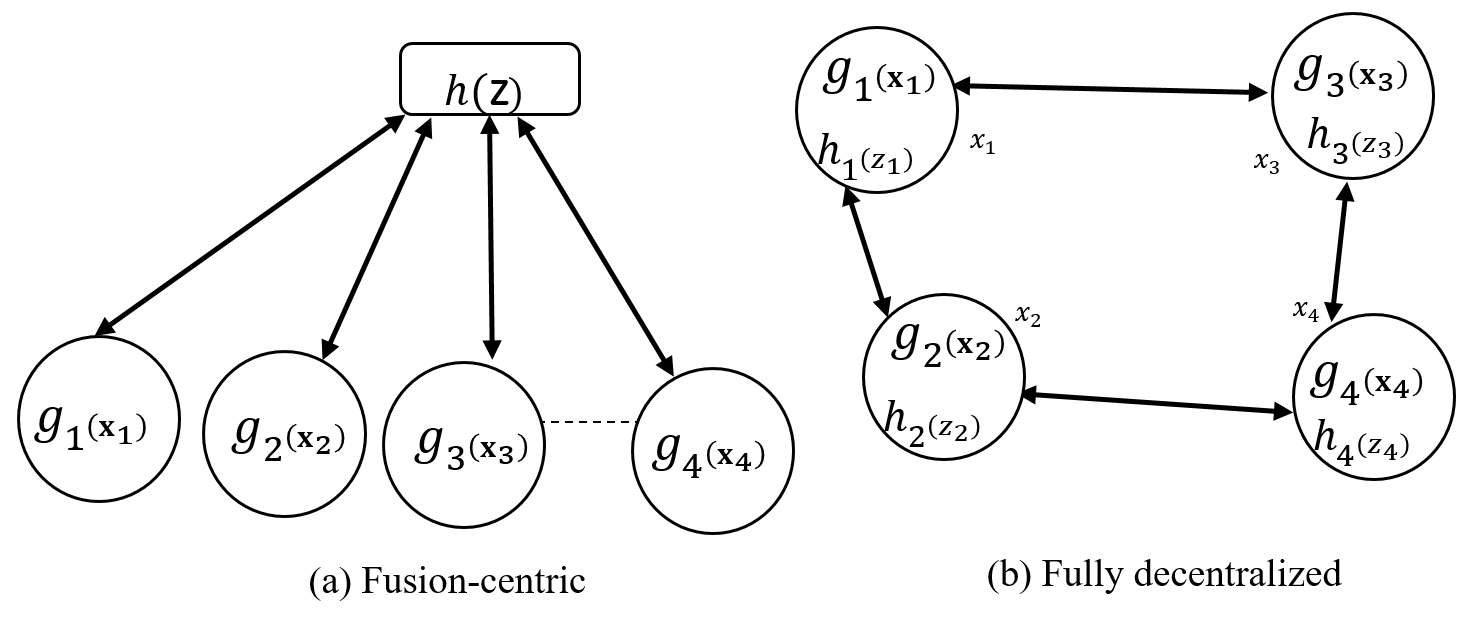}
		\par\end{centering}
	\caption{Distributed optimization architectures}
	\label{admm-architecture-fig}
\end{figure}

The consensus problem in \eqref{unif} will be solved via a successive-convex approximation ADMM algorithm. Of particular interest are distributed implementations where the computational load at individual nodes is minimal. As we shall see later, the formulation in \eqref{unif} is very general and allows for different parallel and distributed ADMM algorithms. The focus here is on two common frameworks: fusion centric and fully decentralized; see Fig. \ref{admm-architecture-fig}.

%\begin{align}
%\mathsf{P} =\underset{\{\x_k\}, \z \in \mathcal{X}}{\min} &\enskip h(\z) + \sum_{k=1}^K g_k(\x_k)\nonumber\\
%\text{s.t. } &\; \x_k = \z \ \ \ \ k = 1, \ldots, K \label{dist-lin:p1c} 
%\end{align}

\subsubsection{Fusion-centric architecture}\label{distributed-def}
The fusion-centric or master-worker architecture has been widely used in data science applications \cite{cevher2014convex,bertsekas1989parallel}, where the aim is to divide the computational load among distributed processors or threads. More generally, such an architecture is useful for any multi-agent system equipped with a central server. In many applications, the loss functions $g_k$ arise from data stored at each node $k$ and the required manipulations, such as the calculation of $\nabla g_k$, constitute the computationally intensive steps that must be carried out at the nodes. In contrast, the fusion center is tasked with the application of the regularization penalty $h$ and overall coordination. 

Associating $g_k$ with the $k$-th node in the network, distributed ADMM methods for such network topologies may be developed by viewing them as special cases of \eqref{unif}. Setting $\pk = \I$, the fusion-centric formulation can be written as

\begin{align} \label{dist-lin:p1c} 
\mathsf{P} =	\underset{\{\x_k\},\z\in \cX}{\text{{minimize}}} &
\begin{array}{c}
\hspace{.2cm}	h(\z) + \sum_{k=1}^K g_k(\x_k) 
\end{array}\\
\text{{subject to}} & \begin{array}[t]{l} 
\hspace{.2cm}\x_k = \z \ \ \ \ k = 1, \ldots, K. 	
\end{array}
\end{align}
Examples of the function $g_k$ include minibatch of loss functions modeling data fidelity, activation function in neural networks such as logit or tanh, and utility functions in network resource allocation. Examples of regularizers include smooth non-convex approximations to $\ell_0$ norm or rank function \cite{gasso2009recovering}, sparsity-promoting convex approximations such as $\ell_1$-norm, nuclear norm, or elastic net regularizer, and indicator functions $\In_{\cX}(.)$ for convex and closed feasible set $\cX$ \cite{oh2013non,friedman2010regularization,hajinezhad2016nestt,bertsekas1989parallel,admm_defn,admm}.

\subsubsection{Fully decentralized architecture}\label{decentralized-def}
Multi-agent systems arising in tactical ad hoc networks, cyber-physical systems, robotic networks, and wireless sensor networks are often spatially distributed and lack a central coordinator or fusion center. In such applications, the estimation, control, or resource allocation tasks must be performed in a fully distributed manner and without reliance on multi-hop communications. Unlike a fusion-centric architecture, fully distributed algorithms are only possible for problems with special structure, such as convexity or partial separability. For instance, consensus-based distributed gradient and subgradient algorithms have been widely applied to convex problems (see e.g. \cite{nedic2009distributed}), while ADMM-based distributed algorithms have been proposed for non-convex partially separable problems \cite{kumar2017asynchronous}. Note that the existing successive convex approximation algorithms either cannot be implemented in a fully distributed manner or have high computation cost due to the need to run nested loops \cite{scutari2018parallel, scutari_part_1, scutari_part_2}.

%[Move this part to intro? ] Multi-agent networked systems arise in a number of engineering disciplines such as tactical ad-hoc networks\cite{martin2005distributed}, environmental monitoring networks \cite{speranzon2006distributed}, multi-robot control and tracking \cite{localization_example, hero}, internet-scale monitoring\cite{mateos2013dynamic,cortes2009distributed}, and large-scale learning \cite{forero2008consensus}. The estimation, resource allocation, and network control tasks required in these applications are often formulated as distributed optimization problems where each node is associated with a local cost function, determined from the set of local and possibly private measurements \cite{admm_defn,duchi2012dual}. Often these applications are constrained by many factors like lack of centralized infrastructure, geographical constraints, resource constraints and societal constraints. Which require decentralized algorithms that do not rely on a fusion center, cluster heads, or multi-hop communication.
%Such framework is widely applicable in network applications such as localization, distributed interference alignment, power system state estimation and high dimensional optimization \cite{kumar2017asynchronous,zheng2018multiagent,assran2018asynchronous,kumar2017distributed}.

We consider a special case of \eqref{unif} where the objective function is partially separable and therefore amenable to a fully decentralized implementation. Specifically, let the network be represented by the undirected graph $\cG = (\cK,\cE)$, where $\cK : = \{1, 2, \ldots, K\}$ denotes the set of agents or nodes, and $\cE$ the set of edges that represent communication links. A node $k \in \cK$ may only communicate with its neighbors $\N'_k := \{j | (j,k) \in \cE\}$. The optimization variable $\x \in \Rn^K$ and each variable $x_k$ is associated with the node $k$. Further, the loss functions are locally coupled, and $g_k(\x_k)$ depends only on the variables $\{x_n | n \in \N_k\}$ where $\N_k := \N'_k \cup \{k\}$. On the other hand, the regularizer and the constraint sets split into nodal components, that is, $h(\z) = \sum_k h_k(z_k)$ and {$\cX:=\times_{k} \cX_k$ where $\cX_k \subset \Rn$}. Introducing variables $\{x_{kj}\}_{j\in\N_k}$ corresponding to every node $k$, the required consensus constraints become $x_{kj} = z_j$. Collecting the variables local to node $k$ into the vector $\x_k \in \Rn^K$ such that $[\x_k]_j = x_{kj}$ for $j \in \N_k$ and zero otherwise, the general form consensus problem becomes 

\begin{align}\label{netw:lin:prob2} 
\mathsf{P} =	\underset{\{\x_k\},\z\in \cX}{\text{{minimize}}} &
\begin{array}{c}
\hspace{.2cm}	h(\z) + \sum_{k=1}^K g_k(\x_k) 
\end{array}\\
\text{{subject to}} & \begin{array}[t]{l}
\hspace{.2cm} x_{kj}=z_j, \ \ \ \ j\in \N_k, k \in \cK.
\end{array}
\end{align}
The formulation \eqref{netw:lin:prob2} is a special case of \eqref{unif} with $K = n$ and $[\P_{k}]_{j\ell} = 1$ for $j = \ell \in \N_k$ and zero otherwise. It can be seen that $\sum_k \P_k$ is invertible and that $g_k$ depends on $\x_k$ through $\P_k\x_k$ only, as required. As an illustration, a four-node network is shown in Fig. \ref{admm-architecture-fig}(b), where each variable $\{x_k\}$ is local to node $k$, and component function $g_k(\x_k)$ at $k$ depends only on variables that are local either to node $k$ or its neighbors. For instance, $x_1$ is local to node $1$ and $g_1(.)$ depends on $\x_1:=[x_1~x_2~x_3~0]^T$. It is remarked that while \eqref{netw:lin:prob2} is written for the special case when $\{x_{kj}\}$ and $\{z_j\}$ are scalar variables, the formulation and the analysis is applicable to the more general case when each node is associated with multiple variables. Accordingly, the notation for generalized inner product and norm will be retained as in \cite{kumar2017asynchronous}.

While the present work focuses on the two distributed architectures, the formulation \eqref{unif} allows other topologies and implementations as well. For instance, \eqref{unif} includes the centralized ADMM formulation where $K = 1$, mixed architecture where local variables are exchanged among nodes while the fusion center takes care of the global variables, and general network topologies where multi-hop neighbors are allowed. 
%and hybrid topologies where nodes are divided into clusters and cluster-heads act as fusion centers for their neighbors while exchanging updates among themselves. 

%%%%%%%%%%%%%%%%%%%%%%%%%%%%%
\section{Distributed Successive Convex Approximation ADMM (DISC-ADMM) }
We begin with detailing the general successive convex approximation ADMM algorithm for solving \eqref{unif}. The distributed implementation details specific to the formulations \eqref{dist-lin:p1c} in \ref{distributed-def} and \eqref{netw:lin:prob2} in \ref{decentralized-def} will be discussed in Sections \ref{distributed-sec} and \ref{decentralized-sec}, respectively.
 Associating dual variables $\{\y_k\}_{k=1}^K$ with the consensus constraints \eqref{unifc}, the augmented Lagrangian can be written as
 
	\begin{align} \label{ual}
	\leg (\{\x_k\}, \z, \{\y_k\}) & = h(\z) + \sum_{k=1}^K {g_k(\x_k)} \nonumber\\
	&\hspace{-1.5cm}+ \sum_{k=1}^K \ip{\pk\y_k}{\x_k - \z} + \frac{\rho}{2}\sum_{k=1}^K \norm{\pk(\z - \x_k)}^2,
	\end{align}
where we have used the fact that $\pk^T = \pk$. Define,	
	\begin{align}
\ell_k(\x_k,\z,\y_k) :=& g_k(\x_k) + \ip{\pk\y_k}{\x_k-\z} \nonumber\\
& + \frac{\rho}{2} \norm{\pk(\x_k - \z)}^2, \label{lk}\\
u(\{\x_k\},\z,\{\y_k\}) :=& h(\z) + \sum_{k=1}^K\ip{\pk\y_k}{\x_k-\z} \nonumber\\
& + \frac{\rho}{2} \sum_{k=1}^K\norm{\pk(\x_k - \z)}^2, \label{u}
\end{align}
so that $\leg(\{\x_k\}, \z, \{\y_k\}) = h(\z) + \sum_{k=1}^K \ell_k(\x_k,\z,\y_k) = u(\x_k,\z,\y_k) + \sum_{k=1}^K g_k(\x_k)$. With these definitions, the classical ADMM algorithm takes the form

	\begin{subequations}\label{unifup}
		\begin{align}
		\z^{t+1} &= \arg \min_{\z \in \Xs} \leg(\{\x_k^t\},\z,\{\y_k^t\}), & \label{zup}\\
		\x_k^{t+1} &= \arg\min_{\x_k} \ell_k(\x_k,\z^{t+1},\y_k^t) & k = 1, \ldots, K,\label{xupd} \\
		\y_k^{t+1} &= \y_k^t + \rho \pk(\x_k^{t+1}-\z^{t+1}) & k = 1, \ldots, K.\label{yup}
		\end{align}
	\end{subequations}

Direct application of the updates in \eqref{unifup} is not viable as both the sub-problems in \eqref{zup} and \eqref{xupd} entail minimization of non-convex functions. Additionally, although the minimization in \eqref{xupd} is separable over the nodes, the required message passing is not straightforward to implement in real-world networks that suffer from the following impediments:
	\begin{itemize}
		\item \textbf{communications delays}: messages transferred between the nodes may get delayed due to wireless impairments or hardware disruptions;
		\item \textbf{computational delays}: in heterogeneous networks, some nodes may be slow intermittently, taking a long time solve the sub-problems in \eqref{zup} and/or \eqref{xupd};
		\item \textbf{asynchronous clocks}: achieving precise synchronization among nodes may not be feasible; 
		\item \textbf{communication losses}: channel impairments or high delays may result in packets being dropped; and 
		\item \textbf{inexact solution}: it may be difficult to solve the sub-problems in \eqref{xupd} and \eqref{zup} exactly.
	\end{itemize}
Towards addressing the computational and communication challenges, we put forth an asynchronous and inexact successive convex approximation ADMM algorithm. The key features of the proposed algorithm include (a) the use of convex surrogate functions in place of the non-convex ones; (b) the ability to handle errors in the update calculations; and (c) the ability to work with delayed iterates allowing asynchronous and parallel implementation. 
	
%%%%%%%%%%%%%%%%%%%%%%%%%%%
\subsection{Simplified updates via surrogate functions}\label{surrogate-def}
Successive convex approximation has been widely applied to a wide range of non-convex problems \cite{mairal2013optimization,lu2018unified,razaviyayn2013unified,mairal2013stochastic,scutari2018parallel}. The main idea in these algorithms is to replace the complicated non-convex functions with simple convex surrogate functions at each iteration. {SCA provides flexibility in tailoring the choice of the surrogate functions to the specific structure of the optimization problem under consideration, and offer a lot of freedom in the algorithmic design\cite{scutari2018parallel,razaviyayn2014parallel}.} In the present case, since both $g_k$ and $h$ may be non-convex, surrogates may be required for both \eqref{xupd} and \eqref{zup}.

Before we proceed further, we write the useful quadratic upper bound property. Let $f$ be any $L$-smooth function function, then for any two points $\x,\tilde{\x}\in \cX $, we have the following quadratic upper bound property
\begin{align}
f({\x}) \leq f(\bar{\x}) + \ip{\nabla f(\bar{\x})}{\x-\bar{\x}} + \frac{L}{2}\norm{\x-\bar{\x}}^2 \label{qs}
\end{align}

For the update in \eqref{zup}, we consider the following class of functions:
\begin{assumption}\label{hsmooth}
	The regularization penalty is expressible as $h = h^s + h^c$ where $h^s$ is a $L_h$-smooth non-convex component while $h^c$ is the convex but possibly non-smooth component.
\end{assumption}
Given $\zb \in \cX$, a convex surrogate function $\h^s(\z;\zb)$ is required for $h^s(\z)$ that adheres to the following assumption:

\begin{assumption}\label{hs}
	For the $L_h$-smooth function $h^s$, the surrogate function $\h^s$ is $L_h$-smooth and convex function, which satisfies
	 \begin{subequations}
		\begin{align}
	 	h(\z) - \h(\z;\zb) &\leq \frac{L_h}{2}\norm{\z-\zb}^2 {+ \omega(\zb)} \tag{\textbf{P1}}\label{h1}\\
		\nabla \h^s(\zb;\zb) &= \nabla h^s(\zb) \tag{\textbf{P2}}\label{h2}\\
		\norm{\nabla \h^s(\z;\zb) - \nabla \h^s(\z';\zb)} &\leq L_h\norm{\z-\z'} \tag{\textbf{P3}}\label{h3}
 		\end{align}
	\end{subequations}
where {$\omega(\zb) := h^s(\zb)-\h^s(\zb;\zb)$ so that the inequality in \eqref{h1} becomes equality at $\z = \zb$}.
\end{assumption}

{Property \eqref{h1} is more flexible than the upper bound requirement in majorization minimization algorithms \cite{sun2017majorization}. For example, as we shall see later, \eqref{h1} allows block-convex functions, not generally handled within the majorization-minimization framework.} Properties \eqref{h2}-\eqref{h3} pertain to the first order behavior of $h^s$ and $\h^s$ only. Henceforth, $\h := \h^s + h^c$ is a convex function and we modify the $\z$-update in \eqref{zup} to utilize $\h(\z;\z^t)$ instead of $h(\z)$.

Each loss function $g_k$ adheres to the following assumption for all $k = 1, 2, \ldots, K$.

\begin{assumption}\label{gsmooth}
	The loss function $g_k$ is $L_g$-smooth.
\end{assumption}

Given $\bar{\x} \in \cX$, the surrogate $\g_k(\x_k; \bar{\x})$ is used for the update in \eqref{xupd} instead of $g_k(\x_k)$. The requirement for the surrogate function is exactly the same as that in Assumption \eqref{hs} and is repeated below with corresponding constants. 
\begin{assumption}\label{ag}
	Given $\x_k,\x_k' \in \Rn^n$ and $\bar{\x}_k\in\cX$, the surrogate $\g_k$ adheres to
	\begin{align}
	g_k(\x_k) - \g_k(\x_k; \bar{\x}) \leq \frac{L_g}{2}\norm{\pk(\x_k-\bar{\x})}^2 {+ \omega_k(\bar{\x})}& \tag{\textbf{P1'}}\label{p1}\\
	\nabla\g_k(\bar{\x}; \bar{\x}) = \nabla g_k(\bar{\x})\hspace{1.3cm}&\tag{\textbf{P2'}}\label{p2} \\
	\norm{\nabla \g_k(\x_k; \bar{\x}) - \nabla \g_k(\x_k'; \bar{\x})} \leq L_g\norm{\pk(\x_k-\x_k')}& \tag{\textbf{P3'}}\label{p3} 
	\end{align} 
where {$\omega_k(\bar{\x}) := g_k(\bar{\x})-\g_k(\bar{\x}_k; \bar{\x})$ so that the inequality in \eqref{p1} becomes an equality at $\x_k = \bar{\x}$. }
\end{assumption}
%	\begin{align}
%g(\x) \leq g(\bar{\x}_k) + \ip{\nabla g(\bar{\x}_k)}{\x-\bar{\x}_k} + \frac{L_g}{2}\norm{\x-\bar{\x}_k}^2 \label{qs}
%\end{align}

Different from \eqref{h1}-\eqref{h3}, both $g_k$ and $\g_k$ depend on $\x_k$ through $\pk\x_k$ so that $\nabla \g_k(\x_k; \bar{\x}) = \pk\nabla \g_k(\x_k; \bar{\x})$. Interestingly, for most loss functions, such surrogate functions are not difficult to find, as the subsequent examples demonstrate. For the examples, the subscript $k$ is dropped for the sake of brevity.

	\begin{exmp} \label{eglip} Consider an $L_g$-smooth function $g$. Then the quadratic upper bound property as in \eqref{qs} suggests that the surrogate function should take the form $\g(\x; \bar{\x}) = g(\bar{\x}) + \ip{\nabla g(\bar{\x})}{\x}$. It can be seen that the surrogate function is linear in $\x$ and satisfies \eqref{h1}-\eqref{h3} with $\omega(\bar{\x}) = -\ip{\nabla g(\bar{\x})}{\bar{\x}}$. Such surrogates have been widely used in big data applications; see e.g. \cite{mairal2013optimization}. 
	\end{exmp}

	\begin{exmp} \label{egdc} Consider a function that can be decomposed as $g(\x)=g^{+}(\x) - g^{-}(\x)$ where $g^+$ and $g^{-}$ are both smooth convex functions. For such functions, the surrogate function $\g(\x; \bar{\x}) = g^{+}(\x) - g^{-}(\bar{\x}) - \ip{\nabla g^{-}(\bar{\x})}{\x}$ satisfies \eqref{h1}-\eqref{h3} with $\omega(\bar{\x}) = \ip{\nabla g^{-}(\bar{\x})}{\bar{\x}}$. %. If the gradients of $g^{+}$ and $g^{-}$ have Lipschitz continuous gradients with constants $L^{+}_g$ and $L^-_g$ respectively, then we have that $L_1 = L_2 = L^{-}_g$ and $L_3 = L^+_g$. 	
	 %\label{egmm} 
	 More generally, let $\g(\x; \bar{\x})$ be a majorizer of $g(\x)$ that satisfies the upper bound and tangent properties, i.e., $\g(\x; \bar{\x}) \geq g(\x)$, $\g(\bar{\x}; \bar{\x}) = g(\bar{\x})$, and $\nabla_{\x}\g(\bar{\x}; \bar{\x}) = \nabla g(\bar{\x})$. Then, $\g(\x; \bar{\x})$ satisfies \eqref{h1}-\eqref{h3} if $\nabla_\x \g(\x; \bar{\x})$ is Lipschitz continuous with respect to $\x$. 
	\end{exmp}
	\begin{exmp} \label{egprod} Consider the function $g(\x) = f_{1}(\x)f_{2}(\x)$ where $f_{1}$ and $f_{2}$ are bounded and smooth convex functions that satisfy 
		\begin{align}
		|f_{i}(\x)| &\leq G \\
		|f_{i}(\x)-f_{i}(\x')| &\leq L\norm{\x-\x'} \\
		\norm{\nabla f_{i}(\x)} &\leq G' \\
		\norm{\nabla f_{i}(\x)- \nabla f_{i}(\x')} &\leq L'\norm{\x-\x'} 
		\end{align}
		for $i = 1, 2$. Then the surrogate $\g(\x; \bar{\x}):=f_{1}(\x)f_{2}(\bar{\x}) + f_{1}(\bar{\x})f_{2}(\x)$ satisfies \eqref{h1}-\eqref{h3} as follows. From the Lipschitz continuity of $f_{1}$ and $f_{2}$, we have that
		\begin{align}
		&(f_{1}(\x)-f_{1}(\bar{\x}))(f_{2}(\x)-f_{2}(\bar{\x})) \leq L^2\norm{\x-\bar{\x}}^2 \\
		\Rightarrow & g(\x)=f_{1}(\x)f_{2}(\x)\leq \g(\x; \bar{\x}) + L^2	\norm{\x-\bar{\x}}^2 - f_1(\bar{\x})f_2(\bar{\x}).\nonumber 
		\end{align}
		which is Property \eqref{h1}. The tangent property \eqref{h2} holds since we have that $\nabla g(\x) = f_{2}(\x)\nabla f_{1}(\x) + f_{1}(\x)\nabla f_{2}(\x) = \nabla_\x (f_{1}(\x)f_{2}(\bar{\x}) + f_{1}(\bar{\x})f_{2}(\x))|_{\bar{\x} = \x}$. Finally, \eqref{h3} may be verified as follows:
		\begin{align}
		&\norm{\nabla_\x \g(\x; \bar{\x}) - \nabla_\x \g(\x'; \bar{\x})} \nonumber\\
		\leq&\norm{(\nabla f_{2}(\x)-\nabla f_{2}(\x')) f_{1}(\bar{\x})} + \norm{(\nabla f_{1}(\x)-\nabla f_{1}(\x')) f_{2}(\bar{\x})} \nonumber\\
		&\leq 2L'G\norm{\x-\x'}
		\end{align} 
		implying that $L_g = \max\{4L'G,2L^2\}$.
	\end{exmp}
\begin{exmp}
	Consider the $L$-smooth function $f(\x_1,\x_2)$ that is block convex in $\x_1$ and $\x_2$ separately but not jointly. The quadratic upper bound for the smooth function $f$ can be written as
	\begin{align}
	&f(\x_1,\x_2) \leq f(\bar{\x}_1; \bar{\x}_2) + \ip{\nabla_{\x_1} f(\bar{\x}_1; \bar{\x}_2)}{\bar{\x}_1-\x_1} \label{fblip}\\
	&+ \ip{\nabla_{\x_2} f(\bar{\x}_1; \bar{\x}_2)}{\bar{\x}_2-\x_2} +\frac{L}{2}\norm{\bar{\x}_1-\x_1}^2 + \frac{L}{2}\norm{\bar{\x}_2-\x_2}^2. \nonumber
	\end{align} 
	From the block-convexity of $f$, we also have that
	\begin{align}
	f(\x_1; \bar{\x}_2) \geq f(\bar{\x}_1; \bar{\x}_2) + \ip{\nabla_{\x_1}f(\bar{\x}_1; \bar{\x}_2)}{\x_1-\bar{\x}_1}\label{fb1}\\
	f(\bar{\x}_1,\x_2) \geq f(\bar{\x}_1; \bar{\x}_2) + \ip{\nabla_{\x_2}f(\bar{\x}_1; \bar{\x}_2)}{\x_2-\bar{\x}_2}\label{fb2}
	\end{align}
	Subtracting \eqref{fb1}-\eqref{fb2} from \eqref{fblip}, we obtain
	\begin{align}
	f(\x_1,\x_2)& -f(\x_1; \bar{\x}_2) - f(\bar{\x}_1,\x_2) \leq -f(\bar{\x}_1; \bar{\x}_2) \nonumber\\
	&+\frac{L}{2}\norm{\bar{\x}_1-\x_1}^2 + \frac{L}{2}\norm{\bar{\x}_2-\x_2}^2
	\end{align}
	suggesting the surrogate $\f(\x_1,\x_2; \bar{\x}_1; \bar{\x}_2) = f(\x_1; \bar{\x}_2) + f(\bar{\x}_1,\x_2)$ with $\omega(\bar{\x}_1; \bar{\x}_2) = -f(\bar{\x}_1; \bar{\x}_2)$. It can also be seen that the gradients of $f$ and $\f$ match for $\x_i = \bar{\x}_i$ for $i = 1, 2$. Since $f$ is smooth, it can be verified that $\f(\x_1,\x_2; \bar{\x}_1; \bar{\x}_2)$ is also smooth with respect to $\x_1$ and $\x_2$. More generally, it can be shown that for the function $f(\{\x_i\}_{i=1}^m)$, the surrogate $\f(\{\x_i\}_{i=1}^m, \{\bar{\x}_i\}_{i=1}^m) = \sum_{i=1}^m f(\x_i, \{\bar{\x}_j\}_{j\neq i})$ satisfies \eqref{h1}-\eqref{h3}. 
\end{exmp}

Summarizing, a large class of surrogate functions adhere to properties \eqref{h1}-\eqref{h3}. In comparison to the surrogate functions \cite{scutari2014parallel} an additional property \eqref{h1} is required. On the other hand, the surrogates $\h$ and $\g_k$ are not required to be strongly convex. Further, the objective function also allows a non-differentiable component $h^c$. Interestingly, surrogates for all the examples considered in \cite{scutari2014parallel} can be constructed so as to adhere to \eqref{h1}-\eqref{h3}.

%Continuing the example for the case when the surrogate function at node $k$ is taken at $\z^{[t]_k}$, observe that for Algorithm 1, only $\{L^i_{1[t]_k}\}_{k=1, i=1}^{K,2}$ are required while for Algorithms 2 and 3, only $\{L^1_{k[t]_k}\}_{k=1,i=1}^{K,2}$ are required. With some abuse of notation, we denote both quantities with $\{L_k^i\}_{i=1}^2$. While such a notation limits the generality of the proof, it helps keep the expressions compact. 

%For the regularizer, we use the surrogate function
%\begin{align}\label{ht}
%\hspace{-2cm}\h(\z,\xb) := \begin{cases} h(\xb) + \ip{\nabla h(\xb)}{\z-\xb} + \frac{L_h}{2}\norm{\z-\xb}_2^2 \\
%&\hspace{-4cm} \text{for $h$ smooth} \\
%h(\z) &\hspace{-4cm} \text{for $h$ convex}
%\end{cases}
%\end{align}
%for all $\z, \xb \in \cX$ 

\subsection{General update rule for DISC-ADMM }\label{sca:update}
An overarching theme of the present work is the design and analysis of the SCA-ADMM algorithm capable of handling errors and delays in the message passing process. Errors are inexorable in distributed settings, arising from communications noise, computational limitations, or from the lack of sufficient data at each node. We consider inexact ADMM updates that utilize an approximate surrogate function $g_k(\x_k; \z_k, \xib_k)$ where $\xib_k \in \Rn^p$ is a random vector that represents the environmental or system state, and $\bar{\x}=\z_k\in \cX$ is the point of majorijation. As a simple example, the approximate surrogate may be constructed as in Example 1 but using an incorrect gradient $\nabla g_k(\z_k) + \e_k$ where $\e_k \in \Rn^n$ is the random error. In order to write the inexact asynchronous DISC-ADMM updates compactly, some new functions are introduced:

	\begin{align}
%	\ell_k(\x_k,\z,\y_k) := &g_k(\x_k) + \ip{\pk\y_k}{\x_k - \z} \nonumber\\
%	& + \frac{\rho}{2} \norm{\pk(\x_k - \z)}^2 \label{lk}\\
	\lt_k(\x_k,\z,\y_k;\zb_k) := &\g_k(\x_k;\zb_k) + \ip{\pk\y_k}{\x_k - \z} \nonumber\\
	& + \frac{\rho}{2} \norm{\pk(\x_k - \z)}^2\label{vk}\\
	\lt_k(\x_k,\z,\y_k;\zb_k,\xib_k):= &\g_k(\x_k;\zb_k,\xib_k) + \ip{\pk\y_k}{\x_k - \z} \nonumber\\
	& + \frac{\rho}{2} \norm{\pk(\x_k - \z)}^2\label{vkt}\\
	\ut(\{\x_k\},\z,\{\y_k\};\zb) :=& \h(\z;\zb) + \sum_{k=1}^K\ip{\pk\y_k}{\x_k-\z} \nonumber\\
	&+ \frac{\rho}{2} \sum_{k=1}^K\norm{\pk(\x_k - \z)}^2. \label{ukt}
	\end{align}
	
Here, $\lt_k$ in \eqref{vk} is the convex surrogate for the local augmented Lagrangian $\ell_k$ in \eqref{lk}. The expression in \eqref{vkt} is an inexact approximation of \eqref{vk} where the random variable $\xib_k$ is the source of the error. 
Likewise, $\ut$ in \eqref{ukt} is the convex surrogate of $u$ in \eqref{u}. With these definitions, the $\z$-update takes the form 
\begin{align}
\z^{t+1} &= \arg\min_{\z\in\cX} \ut(\{\x_k^t\},\z,\{\y_k^t\};\z^t). \label{zupm} 
\end{align}	
That is, the $\z$-update uses the convex surrogate $\h$ calculated at the previous iterate $\z^t$. { More specifically, at the iteration $t+1$, the master node calculates the surrogate function $\tilde{h}(\z;\bar{\z})$ at $\bar{\z}=\z^t$.} The function $\tilde{h}(\z;\bar{\z})$ is evaluated only through $\z^t$ available at the master node, and can be calculated exactly. In the case when $h(\z) = \sum_{k=1}^K h_k(z_k)$, the update splits into $K$ parallel updates to be carried out at each node individually.

For the $\x_k$-update several modifications are introduced so as to make the algorithm resilient to errors, delays, and losses. First, the update utilizes the inexact surrogate $\lt_k$ in \eqref{vkt} in place of $\ell_k$ in \eqref{xupd}. Second, in order to encourage $\x_k^{t+1}$ to stay close to $\x_k^t$, the Bregman divergence term $d_{\phi}(\pk\x_k,\pk\x_k^t)$ is also included in the objective function of \eqref{xupd}. Given any two points $\u$, $\v \in \Rn^n$ the Bregman divergence is given by $d_\phi(\u,\v):=\phi(\u)-\phi(\v)-\inPro{\nabla \phi(\v)}{\u-\v}$ \cite{banerjee2005clustering, wang2015global} where $\phi:\Rn^n \rightarrow \Rn$ is a strictly convex and continuously differentiable function. As a simple example, the choice $\phi(\u) = \frac{1}{2}\norm{\u}^2_2$ translates to $d_\phi(\u,\v) = \frac{1}{2}\norm{\u-\v}_2^2$. Given $\zb \in \cX$, the modified update sub-problem becomes

\begin{align}\label{vkt:breg}
\x_k^{t+1}=\arg\min_{\x_k} \lt_k(\x_k,\z^{t+1}, \y_k^t; \zb,\xib_k) + \frac{1}{\eta}{d_{\phi}(\x_k,\x_k^t)} 
\end{align}	
where the penalty parameter $\eta > 0$. Since $d_\phi(\x_k,\x_k^t)$ is always convex in $\x$, the update sub-problem \eqref{vkt:breg} is convex.

Third, asynchrony is introduced by calculating the surrogate at an older iterate $\z^{[t+1]_k}$ instead of the most recently calculated $\z^{t+1}$. Here, the index $[t+1]_k$ satisfies $t+1-\tau_1 \leq [t+1]_k \leq t+1$ where maximum tolerable delay $\tau_1$ is bounded. For notation ease, we collect all random components realized at time $t$ into the vector $\xib^t$ so that the inexact asynchronous surrogate function used at iteration $t$ is given by $\lt_k(\x_k, \z^{t+1}, \y_k^t; \z^{[t+1]_k}, \xib^{[t+1]_k})$. Fourth, we let the $\x_k$ update to be optional so that only a subset $\S^x_t$ of nodes may perform the update at time $t$. It is remarked that the update must still be carried out often enough and the precise requirements will be detailed later. Summarizing, the $\x_k$ update at node $k$ may be written as

\begin{align}
\x_k^{t+1} = \begin{cases} \arg\min_{\x_k} ~\lt_k(\x_k,\z^{t+1},\y^t_k;\z^{[t+1]_k},\xib^{[t+1]_k}) \\
\hspace{.5cm} + \frac{1}{\eta} d_\phi(\pk\x_k,\pk\x_k^t) \label{xupm} & \hspace{-1.5cm} k \in \S^x_t \\
\x_k^t & \hspace{-1.5cm} k \notin \S^x_t
\end{cases} 
\end{align}

	%Observe that $\lt_k(\x_k,\z,\y_k,\Zb_k,\xib_k)$ is strongly convex in $P_k(\x_k)$ while $\ut(\z,\{\x_k\},\{\y_k\},\x^t_1)$ is strongly convex in $\z$. 
	
Finally, the dual updates are also optional and are given by 
\begin{align}\label{yupm}
\y_k^{t+1} &= \begin{cases} \y_k^t + \rho \pk(\x_k^{t+1}-\z^{t+1}) & k \in \S^x_t\\
\y_k^t & k \notin \S^x_t 
\end{cases}
\end{align}

Before concluding, it is remarked that the gradient error incurred due to the inexact update is given by
\begin{align}
&\e^{t+1}_k:= \nabla \g_k(\x_k^{t+1}; \z^{[t+1]_k},\xib^{[t+1]_k}) - \nabla \g_k(\x_k^{t+1};\z^{[t+1]_k}) \label{ek}\\
&= \nabla \lt_k(\x_k^{t+1},\z^{t+1},\y^t_k;\z^{[t+1]_k},\xib^{[t+1]_k}) \nonumber \\
&\hspace{2.5cm} -\nabla \lt_k(\x_k^{t+1},\z^{t+1},\y^t_k;\z^{[t+1]_k}).\nonumber
\end{align}
The subsequent performance analysis will depend on the average error given by 
\begin{align}
E_T = \frac{1}{T}\sum_{t=1}^T\E\norm{\e_k^t}^2
\end{align}
Different from the stochastic gradient or stochastic ADMM settings, the gradient errors in \eqref{ek} are not required to be independent or unbiased. Consequently, the proposed algorithms are applicable to inexact variants where the error is deterministic, such as the accelerated methods, variance-reduced methods, and algorithms with quantization errors or perturbations. 

\subsection{Connections to existing works}
The proposed DISC-ADMM algorithm bears some resemblance to the ADMM variants in \cite{hajinezhad2016nestt, huang2016stochastic,di2016next,kumar2017asynchronous,sun2017convergence}, {parallel SCA variants in \cite{scutari2018parallel,koppel2018parallel,yang2016parallel,liu2018stochastic}}and \cite{dvurechensky2017gradient}. The novel aspects and key features of the proposed framework in contrast to the existing works are as follows. 
	\begin{itemize}
\item In terms of parallel implementation previous works have analyzed fusion-centric and fully distributed implementations separately; for example, fusion-centric architectures have been analyzed in \cite{hajinezhad2016nestt,admm,zhang2014asynchronous,de2015variance,huo2016asynchronous,chang2016asynchronous,lian2015asynchronous} while fully distributed architectures have been considered in \cite{di2016next,kumar2017asynchronous,yang2016parallel}. On the other hand, the proposed framework unifies the two settings and provides a general performance analysis and application to both. 
\item A number of ADMM algorithms for non-convex problems have been developed to handle independent identically distributed stochastic errors \cite{ghadimi2013stochastic,ghadimi2015accelerated,hajinezhad2016nestt,yang2016parallel}. The inexact framework subsumes these as well as more general error specifications arising in different ADMM variants. 
\item Inexact gradient-based methods for non-convex optimization problems have been recently studied in \cite{sun2017convergence,dvurechensky2017gradient} but the algorithms are centralized. In contrast, developing asynchronous and distributed algorithms is the main focus of the current work. 
\item Existing ADMM variants for solving non-convex problems are either centralized \cite{huang2016stochastic} or distributed but fusion-centric and synchronous \cite{hajinezhad2016nestt}. The proposed algorithm can be implemented in a distributed manner with asynchronous and inexact updates. 
\item The fully distributed architecture for solving non-convex problems in an asynchronous manner has been considered in \cite{kumar2017asynchronous} but using exact gradients and convex regularizers. A parallel algorithm for solving non-convex problems in an asynchronous manner was proposed in \cite{di2016next} and extended to stochastic case \cite{yang2016parallel}, but is restricted to the fusion-centric topology and interdependent, unbiased update errors.

\item The distributed and parallel SCA methods in \cite{scutari2018parallel} are similar in spirit but involve nested loops: the inner loop must run over the network for several iterations before convergence, and results in a single iteration of the outer loop. Additionally, the approach in \cite{scutari2018parallel} is not designed to handle delays and asynchrony, limiting its applicability to real-world distributed systems. Finally, the stochastic variants of \cite{scutari2018parallel}, reported in \cite{yang2016parallel,liu2018stochastic,koppel2018parallel} cannot handle uncertainties arising due to deterministic or 
adversarial sources.
\end{itemize}

In summary, the proposed algorithm with updates in \eqref{zupm}, \eqref{xupm}, and \eqref{yupm} is the first such algorithm for non-convex optimization problems that can be implemented over general network topologies and is capable of tolerating losses, errors, and delays in a seamless fashion. 

Having detailed the general algorithm, the subsequent section discusses various implementation aspects specific to the fusion-centric and fully distributed architectures. While each of the settings may require additional topology-specific assumptions, both algorithms will be analyzed within the unified framework developed here in Sec. \ref{conv}.

\section{Implementation in distributed topologies}
\subsection{Fusion-centric implementation} \label{distributed-sec}
The fusion-centric implementation entails carrying out the updates in \eqref{zupm}, \eqref{xupm}, and \eqref{yupm} with $\pk = \I$. While the algorithm is applicable for general surrogate functions satisfying \eqref{hs} and \eqref{p1}-\eqref{p3}, the updates take the following simple forms for linear surrogate functions (cf. \eqref{qs}):

\begin{subequations}\label{fcup}
\begin{align}
\z^{t+1} &= \px{h^c+\ind_\cX}{\frac{L_h \z^t+\sum_{k=1}^{K} ( \rho\x_k^t + \y_k^t)-\nabla h^s(\z^{t})}{L_h+\rho K}}. \label{dist-lin:prox-z} \\
\x_{k}^{t+1} &=\frac{\rho\z^{t+1} - \y_k^t + \frac{1}{\eta}\x_k^t - \nabla g_k( \z^{[t+1]_k},\xib^{[t+1]_k})}{1/\eta + \rho} \label{dist-lin:x_k update:async} \\
\y_k^{t+1} &= \y_k^t + \rho \left ( \x_k^{t+1} - \z^{t+1} \right ) \label{dist-lin:dualupdate}
\end{align}
\end{subequations}
where the projection operation in \eqref{dist-lin:prox-z} is defined as

\begin{align}
\px{h^c+\ind_\cX}{\u} := \arg\min_{\v \in \cX} h^c(\v) + \frac{K\rho}{2}\norm{\v-\u}^2_2.
\end{align}
and $\nabla g_k( \z^{[t+1]_k},\xib^{[t+1]_k})$ is an inexact gradient of $g_k$ evaluated at the older iterate $\z^{[t+1]_k}$. 

The updates in \eqref{fcup} can be implemented as follows: at some local time $t_k\geq t$, node $k$ receives $\z^t$ from the FC and replies with $\nabla g_k(\z^t, \xib^t)$, possibly after some delay. The FC sends $\z^t$ at each $t$ to all the nodes and waits for a fixed time to receive the individual gradients. At the $t$-th time slot, the FC receives $\nabla g_k(\z^{[t+1]_k},\xib^{[t+1]_k})$ from nodes $k \in \S_t^x$ and applies the updates \eqref{dist-lin:x_k update:async}, \eqref{dist-lin:dualupdate} for these nodes. For all nodes $k \notin \S_t^x$, no updates are calculated so that $x_k^{t+1} = \x_k^t$ and $\y_k^{t+1} = \y_k^t$. For any time slot, if multiple gradients are received, only the most recent one is kept and others are discarded. Likewise, gradients with time stamps older than the most recently available gradient are discarded. Observe that an implication of this policy is that $[t]_k < [t+1]_k$ for all $k \in \S_t$; that is, only `fresh' gradients are used for updates and no reuse is allowed.

The implementation remains exactly the same when other surrogate functions for $h$ satisfying \eqref{hs} are utilized since the $\z$ update occurs at the FC. However, the message passing depends on the surrogate function used for $g_k$. For instance, consider the case when $g_k$ is a difference of convex functions $g_k^{+}$ and $g_k^{-}$ (cf. Example 2) and the inexact surrogate takes the form $\g_k(\x,\z,\xib) = g_k^{+}(\x) - g_k^{-}(\z) - \ip{\nabla g_k(\z,\xib)}{\x-\z}$. For such a surrogate, the $\x_k$-update takes the form

\begin{align}
\x_k^{t+1} = \arg\min_{\x}\ & g_k^{+}(\x) - \ip{\y_k+\nabla g_k^{-}( \z^{[t+1]_k}, \xib^{[t+1]_k})}{\x} \nonumber\\
&+ \frac{\rho}{2}\norm{\x-\z^{t+1}}^2_2
\end{align}
which requires both $\nabla g_k^{-}( \z^{[t+1]_k}, \xib^{[t+1]_k})$ and the full function $g_k^{+}(\cdot)$ to be transmitted from node $k$ to the FC. The proposed fusion-centric inexact asynchronous SCA-ADMM algorithm for linear surrogate functions is summarized in Algorithm \ref{algofc}.

	\begin{algorithm} 
		\caption{Fusion-centric DISC-ADMM}
		\textbf{Updates at the FC}
		\begin{algorithmic}[1]\label{algofc}
			\STATE Initialize $\{ \x_{k}^1,\y_{k}^1 \}$, for all $k=1,\ldots,K$.
			\FOR{$t=1,2,\ldots$}
			\STATE Update $\z^{t+1}$ as in \eqref{dist-lin:prox-z} \label{dist-lin:stoc:conse:async}
			\STATE Transmit $ \z^{t+1}$ to all nodes $k=1,\ldots,K$\\
			\STATE Wait to receive $\{\nabla g_k(\z^{[t+1]_k}, \xib^{[t+1]_k})\}_{k\in\S_t}$ 
			\FOR{$k \in \S_t$}
			\STATE Update $\{\x_k^{t+1}\}$ as in \eqref{dist-lin:x_k update:async} \label{dist-lin:stoc:primal:async}
			\STATE Update $\{\y_{k}^{t+1}\}$ as in \eqref{dist-lin:dualupdate} \label{dist-lin:stoc:dual:async}
			\ENDFOR
			\FOR{$k \notin \S_t$}
			\STATE Set $\x_k^{t+1} = \x_k^t$
			\STATE Set $\y_k^{t+1} = \y_k^t$
			\ENDFOR
			\IF{{$\max_{k}\big\{\norm{\z^{t+1}-\x_k^{t+1}}^2,\norm{\x_k^{t+1}-\x_k^{t}}^2$}\big\}$\leq \epsilon$} \STATE terminate loop \ENDIF
			\ENDFOR\\
			\vspace{2mm}
			\hspace{-5mm}\textbf{Updates at the node $k$}
			\FOR{$t_k = 1$, $\ldots$}
			\STATE (Optional) Receive $\z$ and mark it as $\z^{t_k}$		
			\STATE(Optional) Calculate gradient $ \nabla g_k(\z^{t_k}, \xib^{t_k})$
			\STATE (Optional) Transmit $\nabla\g_k(\z^{t_k}, \xib^{t_k})$ to the FC
			\ENDFOR\\
			\vspace{2mm}
		\end{algorithmic}
	\end{algorithm}

%More generally, the node $k$ receives $\z^t$ from the FC at local time $t_k \geq t$, evaluates $\nabla g_k(\z^t, \xib^{t_k})$, and transmits it back to the FC, possibly incurring delay in the process. Here, the FC waits for a fixed amount of time to receive the inexact gradients from some of the nodes and then proceeds to use the gradient with the most recent time stamp. Denoting the local time at FC by $t$, $\S_t$ be the set of nodes from which the gradient is received. The FC then proceeds to carry out the $\{\x_k,\y_k\}_{k=1}^K$ updates as \eqref{dist-lin:x_k update:async} and \eqref{dist-lin:dualupdate} and proximal consensus update in \eqref{dist-lin:prox-z}. The nodes for which a delayed gradient is received in the current time slot, the update rules in \eqref{dist-lin:x_k update:async} and \eqref{dist-lin:dualupdate} are utilized. For the nodes $k \notin \S_t$, the update rules $\x_k^{t+1}=\x_k^t$ and $\y_k^{t+1}=\y_k^t$ are instead utilized. If multiple gradients are received, the one with the most recent time stamp is utilized, and duplicate gradients, if any, are discarded. Likewise, gradients with time stamp older than the most recent inexact gradient available are discarded. An implication of this policy is that. 

While the detailed performance analysis will be carried out in Section \ref{conv}, the proposed framework requires all delays to be bounded. Within the context of Algorithm \ref{algofc}, such a requirement translates to ensuring that (a) the gradients received at the FC are not too old; and (b) the updates in \eqref{dist-lin:x_k update:async} and \eqref{dist-lin:dualupdate} are applied sufficiently often. The following assumption makes the notion precise. 

\begin{assumption}\label{As-delay}
Recall that $[t+1]_k$ denotes the time stamp of the gradient that is received at time $t$. Likewise, $\S_t$ denotes the set of nodes for which a valid gradient is received and the updates in \eqref{dist-lin:x_k update:async} and \eqref{dist-lin:dualupdate} can be applied. 
\begin{enumerate}
\item[\textbf{a}.] It holds that $t-\tau_1+1\leq [t+1]_k\leq t+1$ if $k\in\S_t$, where the maximum delay $\tau_1$ is bounded.
\item[\textbf{b}.] Let $1\leq j < i \leq T$ be such that $k \in \S_{i}\cap \S_{j}$. Then it must hold that $[i+1]>[j+1]$. In other words, the gradient used for update must always be `fresh'.
\item[\textbf{c}.] Define $(t+1)_k:= \max\{i \leq t+1 \mid k \in \S_{i-1}\}$ as the time of last update, i.e., $\x_k^t=\x_k^{(t)_k}$ and $\y_k^t = \y_k^{(t)_k}$ for all $1\leq t\leq T$. Then it should hold that $t-\tau_2+1 \leq (t+1)_k \leq t+1$, where the maximum delay $\tau_2$ is bounded. 
\end{enumerate}
For brevity, we will simply use $\tau:=\tau_1+\tau_2+1$ as an upper bound for $\tau_1$ and $\tau_2$. 
\end{assumption}
	
\subsection{Decentralized partially separable objective}\label{decentralized-sec}

Recall that for the fully decentralized setting \eqref{netw:lin:prob2}, we have that $K = n$ and $[\pk]_{j\ell} = 1$ for $j=\ell \in \N_k$ and zero otherwise. The separable nature of the problem allows both $\x_k$ and $z_j$ updates to be carried out at individual nodes. Introduce dual variables $\{y_{kj}\}_{j\in\N_k}$ with the constraints in \eqref{netw:lin:prob2} and define $\y_k\in\Rn^n$ such that $[\y_k]_j = y_{kj}$ for $j\in\N_k$ and zero otherwise. Starting with arbitrary $\x_k^1$, $\y_k^1$, and $\z^1$, the $z_j$ update takes the form:

\begin{align}
z_j^{t+1} = \arg \min _{ z_j\in \cX_j } &\h_j(z_j,z_j^t) + \sum _{ k \in \N_j }{ \langle y_{ kj }^t, x_{ kj }^t-z_j \rangle } \nonumber \\ 
& +\frac{\rho}{2}\sum _{ k \in \N_j }{ { (x^t_{ kj }-z_j )^{ 2 } } }\label{netw:lin:z_k_update}
\end{align}

which may again be compactly expressed for linear surrogate functions as in \eqref{dist-lin:prox-z}. Different from the fusion-centric case however, the nodes are allowed to skip the $z_j$ updates and \eqref{netw:lin:z_k_update} may only be applied for nodes $j \in \S_t^z$. For nodes $j \notin \S_t^z$, we simply set $z_j^{t+1} = z_j^t$. The updates for $\x_k$ and $\y_k$ are the same as those in \eqref{xupm} and \eqref{yupm} specialized to the fully distributed problem \eqref{netw:lin:prob2}. As before, for linear surrogates as in Example 1, these updates may be written as
\begin{align}\label{xasynch}
x_{ kj }^{ t+1 }&= \frac{\rho z_j^{t+1}-y_{kj}^t+\tfrac{1}{\eta}x_{kj}^t-[\nabla g_k(\z^{[t+1]_k},\xib^{[t+1]_k})]_j}{\rho+ 1/\eta}\\
y_{ kj }^{ t+1 } &= y^{ t }_{ kj }+\rho\{ x_{ kj }^{ t+1 }-z_{ j }^{ t+1 }\} \label{netw:lin:eq_y_k1}
\end{align}
for $j \in \N_k$ and $k \in \S_t^x$. Here, $[\nabla g_{ k }( {\z_k})]_j :=\frac{\partial}{\partial x_{kj} }g_k(\x_k) \Big|_{\x_k=\z_k}$ for $j \in \N_k$. 
\vspace{.2cm}

In practice, the implementation of the algorithm requires two sub-slots per time slot. At time $t$ at the first sub-slot, nodes $j \in \S_t^z$ evaluate the $z_j$-update as per \eqref{netw:lin:z_k_update} and broadcast $z_j$ to their neighbors. At the end of each sub-slot, a node may or may not have received $z_j$ from their neighbors and may decide to calculate $\nabla g_k(\z^{[t+1]_k},\xib^{[t+1]_k})$. In the second sub-slot, nodes $k \in S_t^x$ update according to \eqref{xasynch} and \eqref{netw:lin:eq_y_k1} while the other nodes simply use $\x_k^{t+1} = \x_k^t$ and $\y_k^{t+1} = \y_k^t$. Subsequently, the node $j$ may decide to broadcast $\{\rho x_{kj}^t+y_{kj}^t \}$ to neighbors $j \in \mathcal{N}_k$. At the next time slot, only nodes that receive updates from all their neighbors carry out the $z_j$-update. The proposed asynchronous algorithm for node $k \in \cK$ is summarized in Algorithm \ref{netw:lin:async_algo}. Various steps designated as `optional' are still constrained by the restrictions in Assumption \ref{As-delay} and the following assumption:
\begin{assumption}\label{fj}
The frequency of updates of each $z_j$ should be bounded away from zero, i.e., $\abs{\{1\leq t\leq T | j\in\S_t^z\}} > Tf_j$ where $f_j \in (0,1]$ is a constant.
\end{assumption}

\begin{algorithm}
	\caption{Fully-decentralized DISC-ADMM}
	\begin{algorithmic}[1]\label{netw:lin:async_algo}
		\STATE Initialize $\{ x_{kj}^1,y_{kj}^1,z_j^1\}$ for all $j \in\N_k$.
		\FOR{$t=1,2,\ldots$}
		\STATE (Optional) Send $\{\rho x_{kj}^t+y_{kj}^t \}$ to neighbors $j \in \mathcal{N}_k$ \label{netw:lin:sendxa1}
		\IF{$\{\rho x_{jk}^t+y_{jk}^t \}$ received from all $j\in\N_k$} \STATE{(Optional) Update $z_{k}^{t+1}$ as in \eqref{netw:lin:z_k_update} and transmit to each $j \in \N_k$} \label{netw:lin:sendza1} \ENDIF
		\IF{$z^{t+1}_j$ not received from some $j \in \N_k$} \STATE{set $z_j^{t+1} = z_j^t$} \ENDIF 
		\STATE (Optional) Calculate gradient $\nabla g_k(\z^{[t+1]_k},\xib^{[t+1]_k})$ \label{opt:grad}
		\FOR{$k\in\S_t^x$}
		\STATE Update $\x_k^{t+1}:=\{x_{kj}^{t+1}\}_{j\in \N_k}$ as in \eqref{xasynch} \label{netw:lin:xupa1}
		\STATE Update the dual variable $y_{kj}^{t+1}$ as in \eqref{netw:lin:eq_y_k1} \label{netw:lin:yupa1}
		\ENDFOR
		\FOR{$k \notin \S_t$}
		\STATE Set $\x_k^{t+1} = \x_k^t$
		\STATE Set $\y_k^{t+1} = \y_k^t$
		\ENDFOR
		\IF{$\norm{\x_k^{t+1}-\x_k^t} \leq \tau_2$ \label{netw:lin:stopca1}}\STATE terminate loop \ENDIF
		\ENDFOR
	\end{algorithmic}
\end{algorithm}

In summary, the salient features of the proposed asynchronous fully-distributed algorithm are as follows.
\begin{itemize}
	\item All resource-intensive steps, such as gradient calculation, proximal function calculation, and transmission of updates are now optional, subject to Assumptions \ref{As-delay} and \ref{fj}. 
	\item Nodes operating on a power budget may only carry out the updates in Steps \ref{netw:lin:xupa1} and \ref{netw:lin:yupa1} at every iteration. When using the old gradient, these updates amount to simple addition/subtraction operations. It is also possible to defer Steps \ref{netw:lin:xupa1} and \ref{netw:lin:yupa1} to a later time slot when the transmission in Step \ref{netw:lin:sendxa1} occurs.
	\item The nodes may implement a timeout mechanism when listening for updates [cf. Steps \ref{netw:lin:sendxa1} and \ref{netw:lin:sendza1}]. The appropriate default actions must be triggered if nothing is heard from one or more neighbors. 
	\item When using the surrogate function in Example 1 and assuming that the gradient calculation requires a fixed number of floating point operations, node $k$ incurs a computational cost of $\mathcal{O}(|\N_k|)$ and a communication cost of $\mathcal{O}(|\N_k|)$ per iteration. 
\end{itemize}

%\col{Note that different from earlier architecture the algorithm outlined in Algorithm \ref{netw:lin:async_algo} does not require any FC and all computations are done locally with variable exchange possible only by communicating through neighboring nodes. Such framework is widely used in applications over wireless sensor networks and communication networks. Often these applications are resource constrained and their infrastructure are ad hoc in nature. The algorithm for such applications should have low computation, communication cost and should be robust to any bounded level of asynchrony. In this direction, note that the Algorithm \ref{netw:lin:async_algo} is fully decentralized where all the update steps \eqref{netw:lin:sendza1} \eqref{netw:lin:xupa1} and \eqref{netw:lin:yupa1} are done locally by each node $k$ per iteration. The algorithm can tolerate any bounded level of asynchrony arising from both computational and communication delay. To reduce the computational and communication burden demanding steps \eqref{opt:grad}, \eqref{netw:lin:sendza1} and \eqref{netw:lin:sendxa1} are made completely optional. The algorithm implementation is cheap as for each iteration the most computationally demanding step is the gradient calculation which translates to a total computational complexity of $\mathcal{O}(|\N_k|)$ for node k per iteration. Also each node require to communicate among neighboring nodes which translates to communication complexity of $\mathcal{O}(|\N_k|)$ for node k per iteration.}

%%%%%%%%%%%%%%%%%%%%%%%%%%%%%%%%%%%%%%%%%
\section{Convergence Analysis}\label{conv}
%%%%%%%%%%%%%%%%%%%%%%%%%%%%%%%%%%%%%%%%%

\subsection{Preliminaries}
We begin with stating some preliminaries that will be useful for the analysis in Sec. \ref{allproofs}.

%\subsubsection{A lower bound for $h^s-\h^s$} Since $h^s$ is $L_h$-smooth from \eqref{hsmooth}, it follows from the quadratic upper bound that for any $\z$, $\zb \in \cX$, 
%\begin{align}\label{h1a}
%h^s(\z) & \geq h^s(\zb) + \ip{\nabla h^s(\z)}{\z-\zb} - \frac{L_h}{2}\norm{\z-\zb}^2_2 
%\end{align}
%Further, since $\h(\z,\zb)$ is also $L_h$-smooth in $\z$, we have that 
%\begin{align}\label{h2a}
%\h^s(\z,\zb) \leq \h^s(\zb,\zb) + \ip{\nabla \h^s(\zb,\zb)}{\z-\zb} + \frac{L_h}{2}\norm{\z-\zb}^2_2
%\end{align}
%Subtracting \eqref{h2a} from \eqref{h1a}, we obtain
%\begin{align}
%h^s(\z) - \h^s(\z,\zb) & \geq \ip{\nabla h^s(\z)-\nabla \h^s(\zb)}{\z-\zb} - L_h\norm{\z-\zb}^2_2 \nonumber\\
%&\geq -2L_h\norm{\z-\zb}^2_2 \label{hsdiff}
%\end{align}

{\subsubsection{A lower bound for $g_k-\g_k$} Since $g_k$ and $\g_k$ are both $L_g$-smooth in $\x_k$, we have that 
\begin{align}\label{g1a}
&g_k(\x_k) \geq g_k(\z_k) + \ip{\nabla g_k(\x_k)}{\x_k-\z_k} - \frac{L_g}{2}\norm{\x_k-\z_k}^2_2\\ 
&\g_k(\x_k;\z_k) \leq \g_k(\z_k;\z_k) + \ip{\nabla \g_k(\z_k,\z_k)}{\x_k-\z_k} \nonumber\\
&\hspace{4cm} + \frac{L_g}{2}\norm{\x_k-\z_k}^2_2\label{g2a}
\end{align}
Subtracting \eqref{g2a} from \eqref{g1a} and using \eqref{p2}, we obtain
\begin{align}
g_k(\x_k) &- \g_k(\x_k;\z_k) \geq g_k(\z_k) - \g_k(\z_k;\z_k) \nonumber\\
&+ \ip{\nabla g_k(\x_k)-\nabla g_k(\z_k)}{\x_k-\z_k} - L_g\norm{\x_k-\z_k}^2_2 \nonumber\\
&\geq \omega(\z_k)-2L_g\norm{\x_k-\z_k}^2_2 \label{gdiff}
\end{align}
where we have used the Cauchy-Schwarz inequality, Assumption \eqref{gsmooth}, and the definition of $\omega_k(\z_k)$.} 

\subsubsection{Optimality condition for \eqref{zupm}} It can be seen that $u$ is strongly convex in $\z$ with parameter at least $\rho\lambda_p$ where $\lambda_p$ is the smallest eigenvalue of $\sum_k \P_k$. Further, the first order optimality condition for \eqref{zupm} can be written as
\begin{align}
\ip{\partial \ut(\z^{t+1},\{\x_k^t\},\{\y_k^t\};\z^t)}{\z-\z^{t+1}} \geq 0 \label{zopt}
\end{align}
for all $\z \in \cX$, where $\partial \h(\z;\zb)=\nabla \h^s(\z;\zb) + \partial h^c(\z)$. Equivalently, defining the normal cone of $\cX$ at $\x$ as ${N}_{\cX}(\x):=\{\v \in \Rn^n \mid \ip{\v}{\z-\x} \leq 0~\forall~\z\in\cX\}$, we can write \eqref{zopt} as
\begin{align}
0 \in \partial \ut(\z^{t+1},\{\x_k^t\},\{\y_k^t\};\z^t) + N_{\cX}(\z^{t+1}). \label{zopt2}
\end{align}

\subsubsection{Optimality condition for \eqref{xupm}} Since $\nabla d_\phi(\u,\v) = \nabla \phi(\u) - \nabla \phi(\v)$, the optimality condition for \eqref{xupm} becomes:
	\begin{multline}\label{xopt}
	\nabla \lt_k(\x_k^{t+1},\z^{t+1},\y^t_k;\z^{[t+1]_k},\xib^{[t+1]_k}) \hspace{1cm} \\
	+ \frac{1}{\eta}\pk\left(\nabla \phi(\pk\x_k^{t+1}) - \nabla \phi(\pk\x_k^t)\right) = 0
	\end{multline}
	Substituting \eqref{vkt} into \eqref{xopt} and using \eqref{yup}, the optimality condition for \eqref{xupm} becomes
	\begin{multline}\label{xopt2}
	\nabla \g_k(\x_k^{t+1};\z^{[t+1]_k}\xib^{[t+1]_k}) + \pk\y_k^{t+1} \hspace{1cm} \\
	+ \frac{1}{\eta}\pk\left(\nabla \phi(\pk\x_k^{t+1}) - \nabla \phi(\pk\x_k^t)\right) = 0
	\end{multline}
which holds for $k\in\S_t^x$ and $t \geq 1$.

\subsubsection{Compact notation for iterate differences} Let us define the consecutive iteration differences as
\begin{align}
\dx_{k}^t &:= \pk(\x_k^{t+1}-\x_k^t) \\
\delta \y_k^t &:= \pk(\y_k^{t+1}-\y_k^t)\\
\dz^t &:= \z^{t+1}-\z^t. 
\end{align}

\subsubsection{Some useful inequalities} Throughout the proof, we will repeatedly use the following inequalities for real numbers $\{a_i\}_{i=1}^m$
\begin{align}
\left(\sum_{i=1}^m a_i\right)^2 \leq m\sum_{i=1}^m a_i^2. \label{cs1}
\end{align}
which is closely related to the Peter-Paul inequality
\begin{align}
2a_1a_2 \leq \varsigma a_1^2 + \frac{a_2^2}{\varsigma} \label{pp}
\end{align}
for any $\varsigma> 0$. 

%The proof begins by establishing bounds on $\norm{\pk(\z^{t+1}-\x_k^{t+1})}$ and $\norm{\pk(\z^{t+1}-\x_k^{t})}$ that will be useful later on. 

Observe that when $k \in \S^x_t$, since $\pk^2 = \pk$, it holds that $ \delta\y_k^t = \pk(\y_k^{t+1}-\y_k^t) = \rho\pk(\x_k^{t+1}-\z^{t+1})$. In the general case however, we have that $\x_k^{t+1} = \x_k^{(t+1)_k}$ so that 
\begin{align}
&\norm{\pk(\z^{t+1}-\x_k^{t+1})}^2 \nonumber\\
&\leq 2\norm{\pk(\z^{t+1}-\z^{(t+1)_k)}}^2 + \frac{2}{\rho^2}\norm{\delta\y_k^{(t+1)_k-1}}^2\nonumber\\
&\leq 2\delta\sum_{\iota = t-\tau_2+1}^t\norm{\dz^\iota}^2 + \frac{2}{\rho^2}\norm{\delta\y_k^{(t+1)_k-1}}^2. \label{zxbound}
\end{align}
Along similar lines it can also be shown that
\begin{align}
&\norm{\pk(\z^{t+1}-\x_k^{t})}^2 \nonumber\\
&\leq 2(\tau_2+1)\sum_{\iota = t-\tau_2}^t\norm{\dz^\iota}^2 + \frac{2}{\rho^2}\norm{\delta\y_k^{(t)_k-1}}^2. \label{zxbound2}
\end{align}

{Finally, for any $L_g$-smooth function $g_k(\x)$ with $L_g$-smooth surrogate $\g_k(\x,\z)$ satisfying properties \eqref{p2}-\eqref{p3}, it holds that
\begin{align}
&\norm{\nabla \g_k(\x,\z) - \nabla \g_k(\x',\z')} = \|\nabla \g_k(\x,\z) - \nabla\g_k(\z,\z) \nonumber\\
&+ \nabla g_k(\z)- \nabla \g_k(\z') +	\nabla \g_k(\z',\z') - \nabla \g_k(\x',\z')\| \nonumber\\
&\leq L_g\norm{\z-\z'} + L_g\norm{\x-\z} + L_g\norm{\x'-\z'} \label{sm2}
\end{align}
and likewise for $h^s$ satisfying \eqref{h2}-\eqref{h3}. Along similar lines, repeated use of triangle inequality can yield bounds of the form 
\begin{align}
&\norm{\nabla \g_k(\x,\z) - \nabla \g_k(\x',\z')} \leq 
L_g\norm{\x-\x'} + 2L_g\norm{\x'-\z''} \nonumber\\
&+L_g\norm{\z-\z''} + L_g\norm{\z'-\z''} + L_g\norm{\z-\z'} \label{sm3}
\end{align}
for any $\z$, $\z'$, $\z'' \in \cX$. }

%\begin{align}
%&\norm{\nabla \g_k(\x,\z) - \nabla \g_k(\x',\z')} = \| \nabla \g_k(\x,\z) - \nabla \g_k(\x',\z) \nonumber\\
%&+ \nabla \g_k(\x',\z) - \nabla \g_k(\z',\z) + \nabla \g_k(\z',\z) - \nabla \g_k(\z,\z) \nonumber\\
%&+ \nabla \g_k(\z) - \nabla \g_k(\z') + \nabla \g_k(\z',\z') - \nabla \g_k(\x',\z')\|\nonumber\\
%&\leq L_g\norm{\x-\x'} + 2L_g\norm{\z-\z'} + 2L_g\norm{\x'-\z'} \label{sm3}
%\end{align}
%and likewise, 
%}

%\begin{align}
%&\norm{\nabla \g_k(\x,\z) - \nabla \g_k(\x',\z')} = \| \nabla \g_k(\x,\z) - \nabla \g_k(\x',\z) \nonumber\\
%&+ \nabla \g_k(\x',\z) - \nabla \g_k(\z'',\z) + \nabla \g_k(\z'',\z) - \nabla \g_k(\z,\z) \nonumber\\
%&+ \nabla \g_k(\z) - \nabla \g_k(\z') + \nabla \g_k(\z',\z') - \nabla \g_k(\z'',\z') \nonumber\\
%& + \nabla \g_k(\z'',\z') - \nabla \g_k(\x',\z')\|\nonumber\\
%&\leq L_g\norm{\x-\x'} + L_g\norm{\z-\z'} + 2L_g\norm{\x'-\z''} \nonumber\\
%&+ L_g\norm{\z-\z''} + L_g\norm{\z'-\z''}\label{sm4}
%\end{align}
%for any $\z$, $\z'$, $z'' \in \cX$. }

\subsection{Main result and proofs} \label{allproofs}
The following theorem summarizes the main result of the paper.

\begin{thm}\label{thm}
The sequence $\{\z^t\}_{t=1}^T$ generated from the updates in \eqref{zupm}, \eqref{xupm}, and \eqref{yupm} is such that for some $1\leq \ell \leq T$, the point $\z^{\ell}$ is $\mathcal{O}(E_T+1/T)$-stationary for constants $T \gg \rho > \max\{18L_g, \frac{L_h+KL_g\tau}{2\lambda_p}\}$ and $T \gg \eta > \varphi$. 
\end{thm}

In other words, the iterates generated by the proposed algorithm are close to stationary if $E_T$ is sufficiently small. As a special case, for the exact asynchronous SCA-ADMM algorithm, the minimum distance from a stationary point decreases as $\mathcal{O}(1/T)$. On the other hand, such a bound is not very useful if the surrogate function error does not diminish and $E_T$ is not sufficiently small. Interestingly however, a number of ADMM variants such as those utilizing variance reduction, averaging, quantization, or perturbation techniques, all result in diminishing errors and can therefore be analyzed within the inexact ADMM framework. 

The proof of Theorem \ref{thm} proceeds along the following lines. We begin with bounding the successive difference of augmented Lagrangian by expressing it as follows
\begin{subequations}\label{legsplit}
	\begin{align} 
	\leg&\left(\{\x_k^{t+1}\}, \z^{t+1}, \{\y_k^{t+1}\} \right) - \leg \left(\{\x_k^t\}, \z^t, \{\y_k^t\} \right) \nonumber \\
	& = \sum_{k=1}^{F} \ell_k(\x_k^{t+1}, \z^{t+1}, \y_k^{t+1}) - \ell_k(\x_k^{t+1},\z^{t+1}, \y_k^t) \label{legsplity} \\
	& \quad + \sum_{k=1}^{F} \ell_k(\x_k^{t+1}, \z^{t+1}, \y_k^t) - \ell_k(\x_k^t, \z^{t+1}, \y_k^t) \label{legsplitx} \\
	& \quad + \leg \left( \{\x_k^t\}, \z^{t+1}, \{\y_k^t\} \right) - \leg \left(\{\x_k^t\}, \z^t, \{\y_k^t\} \right). \label{legsplitz} 
	\end{align}
\end{subequations}
A bound on each of the three terms in \eqref{legsplit} is first developed in Lemmas \ref{lem2} and \ref{lem3}. Subsequently, upper and lower bounds on the augmented Lagrangian are developed in Lemma \ref{lem1}, which serves as the key result required to establish Theorem \ref{thm}. 
For brevity, define 
\vspace{-.2cm}
\begin{align}\label{w1}
B_k^t := \ip{\nabla \g_k(\x_k^{t+1};{\z^{t+1}}) - \nabla \g_k(\x_k^{t+1};{\z^{[t+1]_k},\xib^{[t+1]_k}})}{\dx_k^t} 
\end{align}
for $k\in\S_t^x$ and $B:=\sum_{k\in\S_t^x}B_k^t$. 

The following lemma establishes simple deterministic bounds on the summands in \eqref{legsplit}. 

%%%%%%%%%%%%%%%%%%%%%%%%%%%%%%%%%%%%%%%%%%%%%%%%%%%%%%%%%%%%%%%%%%%%
%%%%%%%%%%%%%%%%%%%%%%%%%%%%%%%%%%%%%%%%%%%%%%%%%%%%%%%%%%%%%%%%%%%%%%%%%%%%
\begin{lem}\label{lem2} The following bounds hold for all $1\leq k\leq K$:
	\begin{subequations}\label{lem2eq}
\begin{align} 
&\sum_{k=1}^K\ell_k \left( \x_k^{t+1}, \z^{t+1}, \y_k^{t+1} \right) - \ell_k \left( \x_k^{t+1}, \z^{t+1}, \y_k^t \right) \label{lem2y}\\
& \hspace{1cm} \;\;= 	\sum_{k=1}^K\frac{1}{\rho} \norm{\delta \y_k^t}^2\nonumber\\
&\ell_k(\x_k^{t+1}, \z^{t+1}, \y_k^t) - \ell_k(\x_k^t, \z^{t+1}, \y_k^t) \label{lem2x}\\
	& \hspace{1cm} \leq B_k^t - \left(\frac{\rho}{2} + \frac{\varphi}{\eta} - 4L_g\right) \norm{\dx_k^t}^2 + \frac{9L_g}{2\rho^2}\norm{\delta\y_k^t}^2 \nonumber\\
&\leg(\{\x_k^t\},\z^{t+1},\{\y_k^t\})-\leg(\{\x_k^t\},\z^t,\{\y_k^t\}) \label{lem2z}\\
	& \hspace{1cm} \leq - \frac{\lambda_p\rho-L_h}{2}\norm{\dz^t}^2 \nonumber
	\end{align}	
	\end{subequations}
\end{lem}
%%%%%%%%%%%%%%%%%%%%%%%%%%%%%%%%%%%%%%%%%%%%%%%%%%%%%%%%%%%%%
%%%%%%%%%%%%%%%%%%%%%%%%%%%%%%%%%%%%%%%%%%%%%%%%%%%%%%%%%%

The following lemma establishes the required bounds for terms $\norm{\y_k^t}^2$ and $B_k^t$ in Lemma \ref{lem2}.

%%%%%%%%%%%%%%%%%%%%%%%%%%%%%%%%%%%%%%%%%%%%%%%%%%%%%%
\begin{lem}\label{lem3} Given $\rho > 10L_g$, the following bounds hold for all $t \geq 1$ and $k\in\S_t^x$:
	\begin{subequations}\label{lem3eq}
		{\begin{align}\label{lem3y}
			&\norm{\delta\y_k^t}^2 \leq 20\norm{\e_k^{t+1}}^2 + 20\norm{\e_k^{(t)_k}}^2 + 60L_g^2\tau\sum_{\iota = t-\tau}^t \norm{\delta\z^\iota}^2\nonumber\\
			&+20\left(\frac{L_{\phi}^2}{\eta^2}+L_g^2\right)\norm{\dx_k^t}^2 +\frac{20L_{\phi}^2}{\eta^2}\norm{\dx_k^{(t)_k-1}}^2 \\
			&B_k^t \leq \frac{\rho+2L_g}{4}\norm{\dx_k^t}^2 + \frac{L_g}{2}\sum_{\iota = t-\tau_1+1}^t\norm{\dz^\iota}^2 + \frac{1}{\rho}\norm{\e_k^{t+1}}^2\label{lem3w}
			\end{align}}
	\end{subequations}
\end{lem}

Having established bounds on each term separately, the following lemma provides upper and lower bounds on $L(\x_k^{T+1}, \z^{T+1}, \y_k^{T+1})$. 

%%%%%%%%%%%%%%%%%%%%%%%%%%%%%%%%%%%%%%%%%%%%%%%%%%%%%%%%%%%%
\begin{lem}\label{lem1}
	{Given $\rho > 10L_g$, the following upper and lower bounds hold for the augmented Lagrangian at time $T+1$:
		\begin{subequations}
			\begin{align}
			\leg&\left(\{\x_k^{T+1}\}, \z^{T+1}, \{\y_k^{T+1}\} \right) \label{lem1a}\\
			&\leq \leg \left(\{\x_k^1\}, \z^1, \{\y_k^1\} \right) + C_e(\rho)\sum_{t=1}^T\sum_{k\in\S^x_t}\norm{\e_k^t}^2 \nonumber\\
			& - C_x(\rho) \sum_{t=1}^T\sum_{k\in\S^x_t} \norm{\dx_k^t}^2 - C_z(\rho)\sum_{t=1}^T \norm{\dz^t}^2 \nonumber\\
			&\Ex{\leg\left(\{\x_k^{T+1}\}, \z^{T+1}, \{\y_k^{T+1}\} \right)} \label{lem1b}\\
			&\geq \mathsf{P} -C_e'(\rho)\sum_{t=1}^T\sum_{k\in\S^x_t}\norm{\e_k^t}^2 \nonumber\\
			& - C_x'(\rho) \sum_{t=1}^T\sum_{k\in\S^x_t} \norm{\dx_k^t}^2 - C_z'(\rho)\sum_{t=1}^T \norm{\dz^t}^2 \nonumber
			\end{align}
		\end{subequations}
		where,
		\begin{align}
		C_x(\rho)&:= \frac{\rho}{4} + \frac{\varphi}{\eta} - \frac{9L_g}{2}-20\upsilon(\rho)\left(\tfrac{2L_\phi^2}{\eta^2} + L_g^2\right)\\
		C_z(\rho)&:= \frac{\lambda_p\rho - L_h-KL_g\tau}{2} -20K\upsilon(\rho)L_g^2\tau^2 \\
		C_e(\rho)&:= 40\upsilon(\rho) + 1/\rho\\
		C_x'(\rho) &:= \frac{5\upsilon'(\rho)}{2}\left(\tfrac{2L_\phi^2}{\eta^2} + L_g^2\right)\\
		C_z'(\rho) &:= 10\upsilon'(\rho)KL_g^2\tau \\
		C_e'(\rho) &:= 5\upsilon'(\rho)
		\end{align}
		with $\upsilon(\rho):=\frac{9L_g}{2\rho^2}+\frac{1}{\rho}$ and $\upsilon'(\rho) = \frac{1}{\rho-3L_g}$.}
\end{lem}
%%%%%%%%%%%%%%%%%%%%%%%%%%%%%%%%%%%%%%%%%
The proofs of Lemmas \ref{lem2}, \ref{lem3}, \ref{lem1}, and Theorem \ref{thm} are carried out in Appendices \ref{proof-lem2}, \ref{proof-lem3}, \ref{proof-lem1}, and \ref{proof-thm}, respectively.
	
%	 Lemma \ref{lem3}, Lemma \ref{proof-lem3}, \ref{proof-lem1} and Lemma \ref{lem1}	are the intermediary results needed for establishing the convergence proof of Theorem \ref{thm}, and their usage for subsequent analysis is self-explanatory. Due to the space limitations, we do not
%include their proofs here. Interested readers are requested to refer the supplementary material for the detailed derivation.}

\section{Conclusion}
{In this paper, we proposed a unified framework for solving distributed non-convex optimization problems under required constraints of big data and networked-system applications. The proposed algorithms are provably convergent and amenable for distributed and decentralized implementation. As a byproduct of our investigation, we also enrich the theory of SCA methods and the inexact gradient methods. We proposed a more general approximation scheme that adds to flexibility in choosing surrogate functions. The presented analysis of integrating generalized inexact gradient framework within the non-convex ADMM holds an independent contribution for the modified gradient literature and has the potential for the extension to other optimization approaches. 
	
	%The asymptotic convergence analysis states that the proposed algorithm satisfies the stationary condition.} 

%the proof for \ref{lem3}, and \ref{lem1} are \ref{proof-lem2}, \ref{proof-lem3} and \ref{proof-lem1}, respectively. Finally, the convergence proof for Theorem \ref{thm} is provided in Appendix \ref{proof-thm}.
%
% Due to space limitations, we do not include their proofs here. Interested readers may refer to the supplementary material of the current paper

\appendices
%%%%%%%%%%%%%%%%%%%%%%%%%%%%%%%%%%%%%%%

\section{Proof of Lemma \ref{lem2} }\label{proof-lem2}

\begin{IEEEproof}[Proof of \eqref{lem2y}] 
	From the dual update \eqref{yupm} we have that
	\begin{align} \label{splity1}
	\ell_k &\left( \x_k^{t+1}, \z^{t+1}, \y_k^{t+1} \right) - \ell_k \left( \x_k^{t+1}, \z^{t+1}, \y_k^t \right) \\
	& = \ip{\delta \y_k^t}{\x_k^{t+1} - \z^{t+1}} = \frac{1}{\rho} \norm{\delta \y_k^t }^2
	\end{align}
	that holds irrespective of whether $k \in\S^x_t$ or not. 
	\end{IEEEproof}

%%%%%%%%%%%%%%%%%%%%%%%%%%%%%%%%%%%%%%%%
\begin{IEEEproof}[Proof of \eqref{lem2x}] Observe that both sides of \eqref{lem2x} are zero for $k \notin \S^x_t$. Therefore, we only need to establish the bound for the case when $k \in \S^x_t$. Using the definitions \eqref{lk} and \eqref{vk} for $\zb_k=\z^{t+1}$, it follows from property \eqref{p1} that
\begin{align}
	\ell_k&(\x_k^{t+1}, \z^{t+1}, \y_k^t) - \lt_k(\x_k^{t+1},\z^{t+1},\y_k^t;\z^{t+1}) \nonumber\\
	&\leq \frac{L_g}{2}\norm{\pk(\x_k^{t+1}-\z^{t+1})}^2 {+ \omega(\z^{t+1})} \label{splitx1}\\
	&\leq \frac{L_g}{2\rho^2}\norm{\delta\y_k^t}^2{ + \omega(\z^{t+1})} \label{splitx1a}
\end{align} 
where \eqref{splitx1a} follows from the update $\delta \y_k^t = \rho\pk(\x_k^{t+1}-\z_k^{t+1})$ for $k \in \S_t^x$. Subtracting $\ell_k(\x_k^t, \z^{t+1}, \y_k^t)$ from both sides, and introducing $\lt_k(\x_k^t,\z^{t+1},\y_k^t;\z^{t+1})$, the inequality becomes
\begin{align}\label{splitx2}
	&\ell_k(\x_k^{t+1}, \z^{t+1}, \y_k^t) - \ell_k(\x_k^t, \z^{t+1}, \y_k^t) \\
	& \leq \lt_k(\x_k^{t+1},\z^{t+1},\y_k^t;\z^{t+1}) - \lt_k(\x_k^t,\z^{t+1},\y_k^t;\z^{t+1}) \nonumber\\
	& + \lt_k(\x_k^t,\z^{t+1},\y_k^t;\z^{t+1}) - \ell_k(\x_k^t, \z^{t+1}, \y_k^t) + \frac{L_g}{2\rho^2}\norm{\delta\y_k^t}^2 \nonumber
\end{align}
Next, observe that $\lt_k(\x_k,\z^{t+1},\y_k^t;\z^{t+1})$ is strongly convex in $\pk\x_k$ or equivalently, the function $\lt_k(\x_k,\z^{t+1},\y_k^t;\z^{t+1}) - \frac{\rho}{2}\norm{\pk\x_k}^2$ is convex. Therefore, the first term on the right of \eqref{splitx2} can be bounded as 
	\begin{align}
	\lt_k&(\x_k^{t+1},\z^{t+1},\y_k^t;{\z^{t+1}}) - \lt_k(\x^t_k,\z^{t+1},\y_k^t;{\z^{t+1}}) \nonumber\\
	&\leq \ip{\nabla \lt_k(\x_k^{t+1},\z^{t+1},\y_k^t;{\z^{t+1}})}{\x_k^{t+1}-\x_k^t} - \frac{\rho}{2} \norm{\dx_k^t}^2 \nonumber\\
	& \hspace{-5mm}= \ip{\nabla \g_k(\x_k^{t+1};{\z^{t+1}}) + \pk\y_k^{t+1}}{\x_k^{t+1}-\x_k^t} - \frac{\rho}{2} \norm{\dx_k^t}^2 
	\end{align}
	Substituting $\pk\y_k^{t+1}$ from the optimality condition \eqref{xopt2} and recalling the definition of $B_k^t$ from \eqref{w1}, we obtain
	\begin{align}
	\lt_k&(\x_k^{t+1},\z^{t+1},\y_k^t;{\z^{t+1}}) - \lt_k(\x^t_k,\z^{t+1},\y_k^t;{\z^{t+1}}) \nonumber\\
	&\leq \ip{\nabla \g_k(\x_k^{t+1};{\z^{t+1}}) - \nabla \g_k(\x_k^{t+1};{\z^{[t+1]_k},\xib^{[t+1]_k}})}{{\dx_k^t}} \nonumber\\
	& - \frac{1}{\eta}\ip{\nabla \phi(\pk\x_k^{t+1}) - \nabla \phi(\pk\x_k^t)}{\dx_k^t} - \frac{\rho}{2} \norm{\dx_k^t}^2 \\
	&\leq B_k^t - \left(\frac{\rho}{2} + \frac{\varphi}{\eta}\right) \norm{\dx_k^t}^2 
	\end{align}
	where the last inequality follows since the function $\phi$ is $\varphi$-convex where $\varphi \geq 0$. For instance, when $\phi(\u) = \frac{1}{2}\norm{\u}_2^2$, we have that $\varphi = 1$. We have also used the fact that $\g_k$ depends on $\x_k$ only through $\pk\x_k$. 
	
For the second term on the right of \eqref{splitx2}, {it follows from \eqref{gdiff} that
\begin{subequations}\label{splitx3}
	\begin{align}
		& \lt_k(\x_k^t,\z^{t+1},\y_k^t;{\z^{t+1}}) - \ell_k(\x_k^t, \z^{t+1}, \y_k^t) \nonumber \\
		& = \g_k(\x_k^t,{\z^{t+1}})- g_k(\x_k^t) \\
		&\leq 2L_g \norm{\pk(\z^{t+1} -\x_k^t)}^2 {- \omega_k({\z^{t+1}})}.\label{splitx3c}\\
		&= 2L_g \norm{\pk(\z^{t+1} -\x_k^{t+1} + \x_k^{t+1} -\x_k^t)}^2 {- \omega_k(\z^{t+1})} \nonumber\\
		&\leq 4L_g \left(\frac{1}{\rho^2}\norm{\delta\y_k^t}^2 + \norm{\dx_k^t}^2\right) {- \omega_k(\z^{t+1})}\label{splitx3d}
	\end{align}
\end{subequations}}
where again the update \eqref{yupm} is used since $k \in \S^x_t$. The bound on \eqref{legsplitx} thus becomes
\begin{align}\label{splitx4}
	&\ell_k(\x_k^{t+1}, \z^{t+1}, \y_k^t) - \ell_k(\x_k^t, \z^{t+1}, \y_k^t) \\
	& \leq B_k^t - \left(\frac{\rho}{2} + \frac{\varphi}{\eta} - 4L_g\right) \norm{\dx_k^t}^2 + \frac{9L_g}{2\rho^2}\norm{\delta\y_k^t}^2. \nonumber
\end{align} 
\end{IEEEproof}
	%%%%%%%%%%%%%%%%%%%%%%%%%%%%%
\begin{IEEEproof}[Proof of \eqref{lem2z}] From the upper bound property in \eqref{h1}, we have that 
\begin{align}\label{splitz1}
	\leg(\{\x_k^t\},\z^{t+1},\{\y_k^t\}) \leq & \ut(\{\x_k^t\},\z^{t+1},\{\y_k^t\};\z^t) \nonumber\\
	& {+ \frac{L_h}{2}\norm{\dz^t}^2 + \omega(\z^t)}
\end{align}
	Subtracting $\leg(\{\x_k^t\},\z^t,\{\y_k^t\}) = \ut(\{\x_k^t\},\z^t,\{\y_k^t\};\z^t)+\omega(\z^t)$ from both sides, the inequality in \eqref{splitz1} becomes
\begin{align}\label{splitz2}
	\leg & (\{\x_k^t\},\z^{t+1},\{\y_k^t\})-\leg(\{\x_k^t\},\z^t,\{\y_k^t\}) \\
	\leq &\ut(\{\x_k^t\},\z^{t+1},\{\y_k^t\};\z^t) - \ut(\{\x_k^t\},\z^t,\{\y_k^t\};\z^t) {+ \frac{L_h}{2}\norm{\dz^t}^2}\nonumber
\end{align}
Since $\z^{t+1}$ minimizes the strongly convex function $\ut(\{\x_k^t\},\z,\{\y_k^t\},\x_1^t)$ over $\z \in \cX$, it follows that
\begin{align}
	\ut&(\{\x_k^t\},\z^{t+1},\{\y_k^t\};\z^t) - \ut(\{\x_k^t\},\z^t,\{\y_k^t\};\z^t) \leq - \frac{\lambda_p\rho}{2}\norm{\dz^t}^2
\nonumber	%\label{splitz3:1}
\end{align} 
	where we have used the fact that $\h$ is convex so that $u$ is $\lambda_p\rho$-convex. Interestingly, the bound remains valid even when some $z_j$ are not updated so that $[\dz^t]_j = 0$ for all $j \notin \S^z_t$. The required bound becomes
\begin{align}
	\leg & (\{\x_k^t\},\z^{t+1},\{\y_k^t\})-\leg(\{\x_k^t\},\z^t,\{\y_k^t\}) \nonumber\\
	&\leq - \frac{\lambda_p\rho-L_h}{2}\norm{\dz^t}^2
\end{align}

\end{IEEEproof}

%%%%%%%%%%%%%%%%%%%%%%%%%%%%%%%%%%%

\section{Proof for Lemma \ref{lem3}}\label{proof-lem3}

\begin{IEEEproof}[Proof of \eqref{lem3y}]
 In general for any $k \in \S^x_t$, it may not necessarily hold that $k \in \S^x_{t-1}$ as well. Indeed, recall that $(t)_k:=\max\{\iota \leq t \mid k \in \S^x_{\iota-1}\}$ and $t-\tau_2 \leq (t)_k \leq t$. Therefore we have that $\y_k^t = \y_k^{(t)_k}$ and $\x_k^t = \x_k^{(t)_k}$ for all $t\geq \tau$. Substituting $\y_k^t$ and $\y_k^{(t)_k}$ from \eqref{xopt2} for $k \in \S^x_t$, we have that
	\begin{align}\label{ybound}
	\delta\y_k^t &= \pk\y_k^{t+1}-\pk\y_k^{(t)_k} \nonumber\\
	= &\nabla\g_k(\x^{(t)_k}_k;\z^{[(t)_k]_k},\xib^{[(t)_k]_k})-\nabla\g_k(\x^{t+1}_k;\z^{[t+1]_k},\xib^{[t+1]_k}) \nonumber\\
	& + \frac{1}{\eta}\pk\left(\nabla \phi(\pk\x_k^{(t)_k}) - \nabla \phi(\pk\x_k^{(t)_k-1})\right) \nonumber\\
	& - \frac{1}{\eta}\pk\left(\nabla \phi(\pk\x_k^{t+1}) - \nabla \phi(\pk\x_k^t)\right)
	\end{align}
	{From \eqref{ek}, the first summand in \eqref{ybound} may be written as $\nabla\g_k(\x^{(t)_k}_k;\z^{[(t)_k]_k})-\nabla\g_k(\x^{t+1}_k;\z^{[t+1]_k}) - \e_k^{t+1} + \e_k^{(t)_k}$. Using \eqref{sm3}, we obtain
	\begin{align}
	&\norm{\nabla\g_k(\x^{(t)_k}_k;\z^{[(t)_k]_k})-\nabla\g_k(\x^{t+1}_k;\z^{[t+1]_k})} \nonumber\\
	&\leq L_g\norm{\dx_k^t}+ 2L_g \norm{\x_k^{t+1}-\z^{t+1}} \nonumber\\
	& + L_g\norm{\z^{t+1}-\z^{[(t)_k]_k}} + L_g\norm{\z^{t+1}-\z^{[t+1]_k}} \nonumber\\
	& + L_g\norm{\z^{[t+1]_k}-\z^{[(t)_k]_k}}
	\end{align}
	where note that $\pk(\x_k^{t+1}-\z^{t+1}) = \frac{1}{\rho}\delta\y_k^t$. Taking squared norm in \eqref{ybound}, and using the inequality in \eqref{cs1} yields
	\begin{align}\label{ybound2}
	&\norm{\delta\y_k^t}^2 \leq \frac{20L_g^2}{\rho^2}\norm{\delta\y_k^t}^2 + 10\norm{\e_k^{t+1}}^2 + 10\norm{\e_k^{(t)_k}}^2 \nonumber\\
	&+10\left(\frac{L_{\phi}^2}{\eta^2}+L_g^2\right)\norm{\dx_k^t}^2 +\frac{10L_{\phi}^2}{\eta^2}\norm{\dx_k^{(t)_k-1}}^2 \nonumber\\
	&+ 10L_g^2\norm{\z^{t+1}-\z^{[(t)_k]_k}}^2 + 10L_g^2\norm{\z^{t+1}-\z^{[t+1]_k}}^2 \nonumber\\
	& + 10L_g^2\norm{\z^{[t+1]_k}-\z^{[(t)_k]_k}}
	\end{align}
	Here for $\rho^2 > 20L_g^2$, the first term can be taken to the right. For the sake of brevity, we assume $\rho > 10L_g$ so that $(1-\tfrac{20L_g^2}{\rho^2})^{-1}<2$, yielding
	\begin{align}\label{ybound3}
	&\norm{\delta\y_k^t}^2 \leq 20\norm{\e_k^{t+1}}^2 + 20\norm{\e_k^{(t)_k}}^2 \nonumber\\
	&+20\left(\frac{L_{\phi}^2}{\eta^2}+L_g^2\right)\norm{\dx_k^t}^2 +\frac{20L_{\phi}^2}{\eta^2}\norm{\dx_k^{(t)_k-1}}^2 \nonumber\\
	&+ 20L_g^2\norm{\z^{t+1}-\z^{[(t)_k]_k}}^2 + 20L_g^2\norm{\z^{t+1}-\z^{[t+1]_k}}^2 \nonumber\\
	& + 20L_g^2\norm{\z^{[t+1]_k}-\z^{[(t)_k]_k}}
	\end{align}}
	It remains to bound the last three terms in \eqref{ybound3}. Recall from Assumption \eqref{As-delay} that $t-\tau_2 \leq (t)_k \leq t$ so that $t-\tau_1-\tau_2 \leq [(t)_k]_k \leq t$ and consequently $t-\tau_2-\tau_1 \leq [t+1]_k-[(t)_k]_k \leq t+1$. Therefore, from \eqref{cs1}, it can be seen that 
	\begin{align}
	\norm{\z^{t+1}-\z^{[(t)_k]_k}}^2 \leq \tau\sum_{\iota = t-\tau}^t \norm{\delta\z^\iota}^2
	\end{align}
	Using similar inequality for the other two terms and introducing additional terms on the right where necessary, we obtain the required bound.
	
	{A caveat here is that \eqref{ybound2} only holds when $(t)_k$ is well-defined. For instance, the index $(1)_k$ is not defined since no update has occurred at or before time $t=1$ for any node $k$. More generally, $(\dx_k^t, \delta \y_k^t)$ are non-zero only for $\S^{-1}_k:=\{t \mid k \in \S_t^x\} = \{t_1, t_2, \ldots\}$ so that $\x_k^{t_1} = \x_k^1$ and $\y_k^{t_1} = \y_k^1$. While \eqref{lem3y} holds for $t\geq t_2$, it does not hold for $t = t_1$ since $\delta\y_k^{t_1} = \pk\y_k^{t_1+1}-\pk\y_k^{1}$. Instead, for $t = t_1$, we have that 
	\begin{align}
	\delta\y_k^{t_1} =& \pk\y_k^{t_1+1}-\pk\y_k^1 \nonumber\\
	& = -\pk\y_k^1 -\nabla\g_k(\x^{t_1+1}_k;\z^{[t_1+1]_k}) - \e_k^{t_1+1} \nonumber\\
	& - \frac{1}{\eta}\pk\left(\nabla \phi(\pk\x_k^{t_1+1}) - \nabla \phi(\pk\x_k^1)\right)\\
	& = \nabla\g_k(\x^1_k,\z^1) -\nabla\g_k(\x^{t_1+1}_k;\z^{[t_1+1]_k}) \nonumber\\
	& - \frac{1}{\eta}\pk\left(\nabla \phi(\pk\x_k^{t_1+1}) - \nabla \phi(\pk\x_k^1)\right) \\
	& - \pk\y_k^1 - \nabla \g_k(\x_k^1,\z^1)- \e_k^{t_1+1}
	\end{align}
	Therefore, it is possible to obtain the same bound as in \eqref{lem3y} if we include additional terms on the right, define $(1)_k = 1$, $\dx_k^{0} = 0$, and $\e_k^1 = \pk\y_k^1 + \nabla \g_k(\x_k^1,\z^1)$.}
	
	\end{IEEEproof}
	
\begin{IEEEproof}[Proof of \eqref{lem3w}]
	{From the definition of the gradient error, we have that
	\begin{align}\label{bbound1}
	B_k^t =& \ip{\nabla \g_k(\x_k^{t+1};\z^{t+1})-\nabla \g_k(\x_k^{t+1};\z^{[t+1]_{k}})}{\dx_k^t} \nonumber\\
	&- \ip{\e_k^{t+1}}{\dx_k^t}
	\end{align}
	The first term in \eqref{bbound1} can be bounded as
	\begin{align}
	&\ip{\nabla \g_k(\x_k^{t+1};\z^{t+1})-\nabla \g_k(\x_k^{t+1};\z^{[t+1]_{k}})}{\dx_k^t} \nonumber\\
	&\leq \norm{\nabla\g_k(\x_k^{t+1};{\z^{t+1}}) - \nabla \g_k(\x_k^{t+1};\z^{[t+1]_k})}\norm{\dx_k^t} \label{gb1}\\
	&\leq L_g\norm{(\z^{t+1}-\z^{[t+1]_k})}\norm{\dx_k^t} \label{gb2}\\
	&\leq L_g\sum_{\iota=t-\tau_1+1}^t \norm{\dz^t}\norm{\dx_k^t} \label{gb3}\\
	&\leq \frac{L_g}{2}\left(\norm{\dx_k^t}^2 + \sum_{\iota=t-\tau_1+1}^t \norm{\dz^t}^2 \right) \label{gradbound}
	\end{align}
	where \eqref{gb1} follows from the Cauchy-Schawarz inequality, \eqref{gb2}	follows from the Lipschitz property \eqref{p3}, \eqref{gb3} follows from the use of triangle inequality and from Assumption \eqref{As-delay}, and \eqref{gradbound} follows from \eqref{pp} with $\varsigma = 1$. Likewise, the second term in \eqref{bbound1} can be bounded from \eqref{pp} with $\varsigma = \rho/2$ yielding
	\begin{align}
	\ip{\e_k^{t+1}}{\dx_k^t} &\leq \frac{1}{\rho}\norm{\e_k^{t+1}}^2 + \frac{\rho}{4}\norm{\dx_k^t}^2 \label{peterpaul}
	\end{align}
	for all $k \in \S^x_t$. Combining \eqref{gradbound} and \eqref{peterpaul}, we obtain the required bound. }
\end{IEEEproof} 

\section{Proof of Lemma \ref{lem1}}\label{proof-lem1}

\begin{IEEEproof}[Proof of \eqref{lem1a}]
Note that in their current forms, the bounds in \eqref{lem3y} and \eqref{lem3w} are only valid for sufficiently large $t$. In order to specify these bounds for all $t$, we assume that the summations involving $\dz^\iota$ always start from $\iota = \max\{1,t-\tau_2-\tau_1\}$. 

{In order to combine the results of Lemmas \ref{lem2} and \ref{lem3}, define and substitute the bounds in \eqref{lem2eq} and \eqref{lem3eq} into \eqref{legsplit} to obtain for all $t\geq 1$:
\begin{align}
	\leg&\left(\{\x_k^{t+1}\}, \z^{t+1}, \{\y_k^{t+1}\} \right) - \leg \left(\{\x_k^t\}, \z^t, \{\y_k^t\} \right) \\
	&\leq \left(20\upsilon(\rho) + \tfrac{1}{\rho}\right)\norm{\e_k^{t+1}}^2 + 20\upsilon(\rho)\norm{\e_k^{(t)_k}}^2 \nonumber\\
	&-\left(\frac{\rho}{4}+\frac{\varphi}{\eta}-\frac{9L_g}{2}-\frac{20\upsilon(\rho)(L_\phi^2+\eta^2L_g^2)}{\eta^2}\right)\sum_{k\in\S_t^x} \norm{\dx_k^t}^2 \nonumber\\
	& + \frac{20L_\phi^2\upsilon(\rho)}{\eta^2}\sum_{k\in\S_t^x}\norm{\dx_k^{(t)_k-1}}^2-\frac{\lambda_p\rho-L_h}{2}\norm{\dz^t}^2 \nonumber \\
	&+20K\upsilon(\rho)L_g^2\tau\sum_{\iota = t-\tau}^t \norm{\dz^\iota}^2 +\frac{KL_g}{2} \sum_{\iota = t-\tau_1+1}^t \norm{\dz^\iota}^2.\nonumber
	\end{align}
Summing over $t = 1, \ldots, T$ and introducing additional terms on the right as necessary, we obtain the required bound. }
\end{IEEEproof}

\begin{IEEEproof}[Proof of \eqref{lem1b}]
	Since $g_k$ is $L_g$-smooth for all $1\leq k\leq K$, we have that
	\begin{align}\label{dist-maj:gkbnd}
	&	g_k(\z^{t+1}) \leq g_k(\x_k^{t+1}) + \ip{\nabla g_k(\x_k^{t+1})}{\pk(\z^{t+1}- \x_k^{t+1})} \nonumber \\ 
	& \qquad + \frac {L_g}{2} \norm{\pk(\x_k^{t+1}-\z^{t+1})}^2\\
	&= g_k(\x_k^{t+1})+\ip{\nabla g_k( \x_k^{t+1})-\nabla	g_k(\z^{t+1})}{\pk(\z^{t+1}- \x_k^{t+1})} \nonumber\\
	& + \ip{\nabla g_k(\z^{t+1})}{\pk(\z^{t+1}- \x_k^{t+1})} + \frac { L_g}{ 2 } \norm{\pk(\x_k^{t+1}-\z^{t+1})}^2 \nonumber\\
	& \leq g_k( \x_k^{t+1})+\ip{\nabla g_k(\z^{t+1})}{\pk(\z^{t+1}- \x_k^{t+1})} \nonumber\\
	&+ \frac{3L_g}{2}\norm{\pk(\x_k^{t+1}-\z^{t+1})}^2. 
	\end{align}
	Therefore, for $t\geq 1$ it holds that
	\begin{align}
	\leg&\left(\{\x_k^{t+1}\}, \z^{t+1}, \{\y_k^{t+1}\} \right) \geq h(\z^{t+1})+\sum _{k=1}^K g_k(\z^{t+1}) \nonumber\\
	&+ \sum_{k=1}^K\ip{\y_k^{t+1}+\nabla g_k(\z^{t+1})}{\pk(\x_k^{t+1}-\z^{t+1})} \nonumber\\
	& + \frac {\rho-3 L_g}{2} \sum_{k=1}^K\norm{\pk(\x_k^{t+1}-\z^{t+1})}^{2}. \label{dist-maj:lbounda}
	\end{align}
	Since $\z^{t+1}\in\cX$, we have that $h(\z^{t+1}) + \sum_k g_k(\z^{t+1}) \geq \mathsf{P}$. Since $\rho' := \rho - 3L_g \geq 0$, the Peter-Paul inequality implies that
	\begin{align}
	&\ip{\y_k^{t+1}+\nabla g_k(\z^{t+1})}{\pk(\x_k^{t+1}-\z^{t+1})} \\
	&\hspace{5mm}\geq -\frac{1}{2\rho'}\norm{\y_k^{t+1}+\nabla g_k(\z^{t+1})} - \frac{\rho'}{2}\norm{\pk(\x_k^{t+1}-\z^{t+1})}\nonumber
	\end{align}
	which allows us to cancel the last term in \eqref{dist-maj:lbounda}, and yields
	\begin{align}
	\leg&\left(\{\x_k^{t+1}\}, \z^{t+1}, \{\y_k^{t+1}\} \right) \nonumber\\
	&\geq \mathsf{P} - \frac{1}{2\rho'}\sum_{k=1}^K\norm{\pk\y_k^{t+1}+\nabla g_k(\z^{t+1})}^2 \label{lbound0}
	\end{align}
	
	{It remains to develop an upper bound on $\norm{\pk\y_k^{t+1}+\nabla g_k(\z^{t+1})}^2$ for each $1\leq k\leq K$. For any $k$, recall that $\y_k^{t+1} = \y_k^{(t+1)_k}$ regardless of whether $k \in S_{t}^x$ or not. Using \eqref{xopt2} for $t$ sufficiently large so that $(t+1)_k$ is well-defined, we obtain
	\begin{align}
	\pk\y_k^{t+1} &= \pk\y_k^{({t+1})_k} = -\nabla \g_k(\x_k^{({t+1})_k};\z^{[({t+1})_k]_k}) - \e_k^{({t+1})_k} \nonumber\\
	&\hspace{-5mm} - \frac{1}{\eta}\pk\left(\nabla \phi(\pk\x_k^{({t+1})_k}) - \nabla \phi(\pk\x_k^{({t+1})_k-1})\right).
	\end{align}
	Adding $\nabla g_k(\z^{t+1})$ on both sides, taking norm, and using the $L_g$-smoothness property of $g_k$ and $\g_k$ (Properties \eqref{p2} and \eqref{p3}, respectively), we obtain
	\begin{align} \label{lbound1}
	&\norm{\pk\y_k^{t+1}+\nabla g_k(\z^{t+1})} \leq \frac{L_\phi}{\eta}\norm{\dx_k^{({t+1})_k-1}} + \norm{\e_k^{({t+1})_k}} \nonumber\\
	&+L_g\norm{\z^{t+1}-\z^{[({t+1})_k]_k}} + L_g\norm{\pk\x_k^{(t+1)_k}-\pk\z^{[({t+1})_k]_k}}
	\end{align}
	Here, the last term can be bounded using triangle inequality as
	\begin{align}
	 &\norm{\pk\x_k^{(t+1)_k} - \pk\z^{[(t+1)_k]_k}} \nonumber\\
	 &\leq \frac{1}{\rho}\norm{\delta\y_k^{(t+1)_k}} +\norm{\z_k^{(t+1)_k} - \z^{[(t+1)_k]_k}} \label{lbound2}
	\end{align}
	where we have used the update rule \eqref{yup}. Combining the results in \eqref{lbound1} and \eqref{lbound2}, observing that there are a total of five terms on right, and using \eqref{cs1}, we obtain
	\begin{align}
	&\norm{\pk\y_k^{t+1}+\nabla g_k(\z^{t+1})}^2 \leq \frac{5L_\phi^2}{\eta^2}\norm{\dx_k^{(t+1)_k}}^2 + 5\norm{\e_k^{(t)_k}}^2 \nonumber\\
	&+10L_g^2\tau\sum_{\iota=t-\tau}^{t}\norm{\dz^\iota}^2+ \frac{5L_g^2}{\rho^2} \norm{\delta\y_k^{(t+1)_k}}^2 \label{ypg}
	\end{align}
	While the right-hand side contains only a few terms, additional summands can be included so that $t$ ranges from $1$ to $T$ as follows:	
	\begin{align}
	&\norm{\pk\y_k^{t+1}+\nabla g_k(\z^{t+1})}^2 \leq 10L_g^2\tau\sum_{t=1}^T\norm{\dz^t}^2\nonumber\\
	&+ 5\sum_{t=1}^T\left(\frac{L_\phi^2}{\eta^2}\norm{\dx_k^t}^2 + \frac{L_g^2}{\rho^2} \norm{\delta\y_k^t}^2 + \norm{\e_k^t}^2\right) \\
	&\leq 5\left(\frac{L_\phi^2}{\eta^2}(1+\frac{40L_g^2}{\rho^2}) + \frac{20L_g^4}{\rho^2}\right)\sum_{t=1}^T\norm{\dx_k^t}^2 \nonumber\\
	&+10L_g^2\tau\left(1+\frac{30L_g^2}{\rho^2}\right)\sum_{t=1}^T\norm{\dz^t}^2 \nonumber\\
	& + 5\left(1+\frac{40L_g^2}{\rho^2}\right) \sum_{t\in \S^{-1}_k}\norm{\e_k^t}^2 
	\end{align}
	where we have used the result from \eqref{lem3y}. For brevity, since $\rho > 10L_g$, we simply use the bounds $30L_g^2/\rho^2 < 40L_g^2/\rho^2 < 1$ and sum over all $1\leq k \leq K$ to obtain
	\begin{align}
	&\norm{\pk\y_k^{t+1}+\nabla g_k(\z^{t+1})}^2 \leq 10\sum_{t=1}^T\sum_{k\in\S_t^x} \norm{\e_k^t}^2 \label{pkybound}\\
	&+ 20 KL_g^2\tau\sum_{t=1}^T\norm{\z^t}^2 + 5\left(\frac{2L_\phi^2}{\eta^2} + L_g^2\right)\sum_{t=1}^T\sum_{k\in\S_t^x}\norm{\dx_k^t}^2\nonumber
	\end{align}
	The required bound is obtained from plugging \eqref{pkybound} into \eqref{lbound0}. }
\end{IEEEproof}

\section{Proof of Theorem \ref{thm}}\label{proof-thm} 

Lemma \ref{lem1} provides the key ingredients required to prove Theorem 1. The proof of Theorem \ref{thm} comprises of two key steps. First bounds developed in Lemma \ref{lem1} are directly used to obtain the bounds on consecutive Lagrangian differences. Second, a bound on the expected value of the minimum subgradient norm is developed, again using various results in Lemma \ref{lem1}.
	
\begin{IEEEproof}
Letting $\leg_1:=\leg(\{\x_k^{1}\}, \z^{1}, \{\y_k^{1}\})$, taking expectations in \eqref{lem1a}-\eqref{lem1b}, and combining, we obtain
\begin{align}
	\mathsf{P} &\leq \leg_1 + (C_e(\rho)+C_e'(\rho))T E_T \nonumber\\
	& - (C_x(\rho)-C_x'(\rho))\sum_{t=1}^T\sum_{k\in\S^x_t}\E\norm{\dx_k^t}^2 \nonumber\\
	& - (C_z(\rho)-C_z'(\rho))\sum_{t=1}^T\E\norm{\dz^t}^2 \label{dist-lin:pl1}
\end{align} 
Next, observe that it is always possible to choose $\rho$ and $\eta$ sufficiently large such that both $(C_x(\rho)-C_x'(\rho))$ and $(C_z(\rho)-C_z'(\rho))$ are positive. This is because for $\rho$ sufficiently large both $\upsilon(\rho)$ and $\upsilon'(\rho)$ go to zero. In general, suitable values of $\rho_{\min}$ and $\eta_{\min}$ may be chosen such that these terms are positive. {For the sake of brevity, we specify $\eta = \varphi$ and assume that $\rho > \rho_{\min}$ where
\begin{align}
\rho_{\min} = \max\{18L_g, \frac{L_h+KL_g\tau}{2\lambda_p}\}.
\end{align}
which is obtained by ignoring the terms depending on $\frac{1}{\rho}$ or $\frac{1}{\rho-3L_g}$. Having chosen $\rho$ and $\eta$ appropriately, denote constants $C_e:=C_e(\rho)+C_e'(\rho)$, $C_x:=C_x(\rho)-C_x'(\rho)$, and $C_z := C_z(\rho)-C_z'(\rho)$. } Rearranging \eqref{dist-lin:pl1}, we obtain
	\begin{align}
\frac{1}{T}\sum_{t=1}^T &\left(C_x\sum_{k\in\S^x_t}\E\norm{\dx_k^t}^2 + C_z\E\norm{\dz^t}^2\right) \nonumber\\
& \leq C_e E_T + \frac{\leg_1 - \mathsf{P}}{T}.
	\end{align}
It follows therefore that
	\begin{align}
	\frac{1}{T}\sum_{t=1}^T \sum_{k\in\S_t^x}\E\norm{\x_k^{t+1}-\x_k^t}^2 \leq \ot\label{xkrate}\\
	\frac{1}{T}\sum_{t=1}^T\sum_{j\in\S_t^z} \E\norm{z_j^{t+1}-z_j^t}^2 \leq \ot \label{zrate}
	\end{align}
Since $t-(t)_k \leq \tau$, it holds that $\abs{\S^{-1}_k} = \mathcal{O}(T)$ for all $1\leq k \leq K$. Consequently, there exists $\ell \in \S^{-1}_k$ such that $\E\norm{\x_k^{\ell+1}-\x_k^\ell}^2 \leq \ot$ for all $1\leq k\leq K$. Likewise, from Assumption \ref{fj}, it follows that there exists $\ell\in\S_t^z$ such that $\E\norm{z_j^{\ell+1}-z_j^\ell}^2 \leq \ot$. 
	
Recalling that $\h = h^c + \h^s$, it follows from \eqref{zopt2} that there exists $\v \in \partial h^c(\z^{t+1}) + N_{\cX}(\z^{t+1})$ such that 
\begin{align}
&\v = -\nabla \h^s(\z^{t+1},\z^t) + \sum_{k=1}^K \pk\y_k^t + \rho\sum_{k=1}^K \pk(\x_k^{t}-\z^{t+1}) \nonumber\\
& = - \nabla\h^s(\z^{t+1},\z^t) +\sum_{k=1}^K \pk\y_k^{t+1} + \rho\sum_{k=1}^K \pk(\x_k^{t}-\x_k^{t+1}).\label{vdef}
\end{align}
Therefore we have that
\begin{align}
&\vartheta_{t+1}\nonumber\\
&:=\hspace{-3mm}\min_{\v \in \partial h^c(\z^{t+1}) + N_{\cX}(\z^{t+1})} \norm{\v + \nabla h^s(\z^{t+1}) + \sum_{k=1}^K {\nabla} g_k(\z^{t+1})}^2 \nonumber\\
&\leq 3\norm{\nabla h^s(\z^{t+1}) - \nabla\h^s(\z^{t+1},\z^t)}^2 \label{varthbound}\\
& +3K\sum_{k=1}^K \norm{\pk\y_k^{t+1} + \nabla g_k(\z^{t+1})}^2 + 3\rho^2 K\sum_{k=1}^K \norm{\dx_k^t}^2 \nonumber
\end{align}
{Recall from \eqref{h2}-\eqref{h3} that $\norm{\nabla h^s(\z^{t+1}) - \nabla\h^s(\z^{t+1},\z^t)} \leq 2L_h\norm{\dz^t}$. Further, as in \eqref{lbound2}, we have that
\begin{align}
&\sum_{t=1}^T\sum_{k=1}^K\norm{\pk\y_k^{t+1} + \nabla g_k(\z^{t+1})}^2 \leq 5\sum_{t=1}^T\sum_{k\in\S_t^x}\norm{\e_k^t}^2 \nonumber\\
&+\frac{5L_\phi^2}{\eta^2}\sum_{t=1}^T\sum_{k\in\S_t^x}\norm{\dx_k^t}^2+10KL_g^2\tau^2\sum_{t=1}^{T}\norm{\dz^t}^2 \nonumber\\
&+ \frac{5L_g^2}{\rho^2} \sum_{t=1}^T\sum_{k\in\S_t^x}\norm{\delta\y_k^t}^2 \\
& \leq \left(\frac{7L_\phi^2}{\eta^2} + L_g^2\right)\sum_{t=1}^T\sum_{k\in\S_t^x}\norm{\dx_k^t}^2 \nonumber\\
&+11KL_g^2\tau^2\sum_{t=1}^{T}\norm{\dz^t}^2+6\sum_{t=1}^T\sum_{k\in\S_t^x}\norm{\e_k^t}^2 \label{pkybound2}
\end{align}
where we have used \eqref{lem3y} and the fact that $\rho > 18L_g$. From \eqref{varthbound} and \eqref{pkybound2}, we obtain, 
\begin{align}
&\sum_{t=1}^T\vartheta_t \leq \left(33K^2L_g^2\tau^2+12L_h^2\right)\sum_{t=1}^T\norm{\dz^t}^2 \nonumber\\
&+ 3K\left(\frac{7L_\phi^2}{\eta^2} + L_g^2 + \rho^2\right)\sum_{t=1}^T\sum_{k\in\S_t^x}\norm{\dx_k^t}^2 \nonumber\\
&+18K\sum_{t=1}^T\sum_{k\in\S_t^x}\norm{\e_k^t}^2
\end{align}}
Dividing both sides by $T$, taking expectation, and using \eqref{xkrate}-\eqref{zrate}, it can be seen that 
\begin{align}
\frac{1}{T}\sum_{t=1}^T \Ex{\vartheta_t} \leq \mathcal{O}\left(E_T+\frac{1}{T}\right)
\end{align}
Equivalently, there exists $1\leq \ell \leq T$ such that $\z^\ell$ is an $\epsilon$-stationary point of the original problem where $\epsilon = \mathcal{O}(E_T+1/T)$. 
\end{IEEEproof}

%
%In this paper, we proposed a unified algorithmic framework for solving a genric non-convex optimization problem arising in novel algorithmic framework solving a generic non-convex smooth objective function and non-smooth convex or smooth non-convex is considered. A more flexible algorithm is developed by considering very general parallel optimization formulation. The parallel formulation considers widely used popular architectures for machine learning and multi-agent network applications. The developed algorithms are parallel, asynchronous and robust to delay and uncertainties. 
%
%
%
%The notion of inexact gradients 

%%%%%%%%%%%%%%%%%%%%%%%%%%%%%%%%%%%%%%%%%%%%%%%%%%%%%
\bibliographystyle{IEEEtran}
\bibliography{IEEEabrv,refs}
%%%%%%%%%%%%%%%%%%%%%%%%%%%%%%%%%%%%%%%%%%%%%%%%%%%%%%

\end{document}